\documentclass{elsarticle}
\usepackage[a4paper,total={6in,9in}]{geometry}
\usepackage{lipsum}
\usepackage{graphicx}
\usepackage{amssymb}
\usepackage{color}
\usepackage{amsmath}
\makeatletter
\def\ps@pprintTitle{%
 \let\@oddhead\@empty
 \let\@evenhead\@empty
 \def\@oddfoot{}%
 \let\@evenfoot\@oddfoot}
\makeatother
\usepackage{url}

\usepackage{lineno}

\begin{document}

\title{The Impact of Temperature and Isolation on COVID-19 in India : A Mathematical Modelling approach}

\author[sai]{D Bhanu Prakash}
\author[sai]{Bishal Chhetri}

\author[sai]{D K K Vamsi \corref{mycorrespondingauthor}}
\cortext[mycorrespondingauthor]{Corresponding author}
\ead{dkkvamsi@sssihl.edu.in}

\author[sai]{Balasubramanian S}
\author[sai,sweden]{Carani B Sanjeevi\fnref{myfootnote}}
\fntext[myfootnote]{Vice-Chancellor, Sri Sathya Sai Institute of Higher Learning -  SSSIHL, India  \\  Email Address : sanjeevi.carani@sssihl.edu.in, sanjeevi.carani@ki.se}

\address[sai]{Department of Mathematics and Computer Science, Sri Sathya Sai Institute of Higher Learning - SSSIHL, India}
\address[sweden]{Department of Medicine, Karolinska Institute, Stockholm, Sweden }

\date{\today}

\begin{abstract}
The dynamics of COVID-19 in India are captured using a set of delay differential equations by dividing a constant population into six compartments. The equilibrium points are calculated and stability analysis is performed. Sensitivity analysis is performed on the parameters of the model. Bifurcation analysis is performed and the critical delay is calculated. By formulating the spread parameter as a function of temperature, the impact of temperature on the population is studied. We concluded that with the decrease in temperature, the average infections in the population increases. In view of the coming winter season in India, there will be an increase in new infections. This model falls in line with the characteristics that increase in isolation delay increases average infections in the population.
\end{abstract}


\maketitle

\section{Introduction}

The pandemic of COVID-19 has spread its tentacles across the world. Large number of countries are seeing a spike in the new infections and a few experiencing second wave of the pandemic. Moving into the winter season where influenza spreads rapidly and vaccine not being available to everyone in the near future, the impact of temperature and isolation is an important study in fighting COVID-19. \\

Various models were framed to understand the dynamics in population scale in India\cite{main,main1,main2,main3}. Since COVID-19 is asymptomatic in many cases, immediate isolation of new infections is not practical. The delay involved in identification and effective isolation plays a major role in flattening the virus curve. A few models incorporated the impact delay on COVID-19 in India\cite{delay,delay1,delay2,delay3}. Since the SARS-CoV-2 virus, responsible for COVID-19, spreads faster through air, the temperature in environment might play a role in spreading the virus. The effect of temperature on COVID-19 is studied in a few papers\cite{temp1,temp2,temp3,temp4}.\\

In this paper, we tried to capture the dynamics of COVID-19 in India using a system of delay differential equations from the inspiration given in \cite{nature}. The effect of temperature is also incorporated using the spread parameters. The rest of the paper goes in this fashion. In Section \ref{sec2}, we proposed the model along with the detailed explanations of delay, temperature and other parameters. In Section \ref{sec3}, a few important indicators like reproduction number and equilibrium points has been calculated. The stability of equilibrium points is discussed along with the Hopf bifurcation analysis in Section \ref{sec4}. In Section \ref{sec8}, the sensitivity analysis of sensitive parameters is given. The sensitivity analysis for the rest of the parameters is given in Appendix 1. Section \ref{sec5} deals with numerical simulations. Section \ref{sec6} presents the impact of temperature and isolation on COVID-19 in India. In Section \ref{sec7}, conclusions have been given. \\

\section{Mathematical Model}
\label{sec2}

We frame the model equations by considering the fraction of the people in each category. We assume a constant population(N) and divide it into six compartments as :

\begin{itemize}
\item S - fraction of the total population that is healthy and has never caught the infection.
\item E - fraction of the total population that is exposed to infection.
\item I - fraction of the total population infected by the virus and undetected.
\item Q - fraction of the total population that are found positive in the test and either hospitalized or quarantined.
\item R - fraction of the total population that has recovered from the infection.
\item D - fraction of the total population that are extinct due to the infection.
\end{itemize}

The proposed model is :

\begin{align*}
\frac{\mathrm{d} S(t)}{\mathrm{d} t} &= \Omega -\beta (T) S(t) I(t) - \mu S(t)  \\
\frac{\mathrm{d} E(t)}{\mathrm{d} t} &= \beta (T) S(t) I(t) - \epsilon E(t) - \mu E(t)  \\
\frac{\mathrm{d} I(t)} {\mathrm{d} t} &= \epsilon E(t) - \gamma I(t) - p e^{-\gamma \tau} I(t - \tau) - \mu I(t)  \hspace{4cm}(2.1)\\
\frac{\mathrm{d} Q(t)}{\mathrm{d} t} &= p e^{-\gamma \tau} I(t - \tau) - \rho(1-\alpha) Q(t-\kappa) - \delta \alpha Q(t) - \mu Q(t) \\
\frac{\mathrm{d} R(t)}{\mathrm{d} t} &= \gamma I(t) + \rho (1-\alpha) Q(t-\kappa) - \mu R(t) \\
\frac{\mathrm{d} D(t)}{\mathrm{d} t} &= \delta \alpha Q(t) - \mu D(t) \label{eqn1} 
\end{align*}

Infectious hosts infect their neighbors at rate $\beta$. $\beta$ is the expected amount of people an infected person infects per day. That infected person can infect only S(t) people. Hence the change of S(t) to the next day = -$\beta$ S(t) I(t). Exposed individuals become infected at the rate $\epsilon$. If an infected individual remains infectious for $\tau$ units of time without having recovered, it enters a new state, Q (for quarantine), with the detection probability p. Among I(t-$\tau$) infected persons at time t-$\tau$, only a fraction of them ($e^{-\gamma \tau}$) remains infectious and they are moved to quarantine Q at the rate p. After staying in quarantine for $\kappa$ days, individuals in Q compartment move to R at the rate $\rho$ (1-$\alpha$) , where $\rho$ is the rate at which people recover and $\alpha$ is the death rate. Individuals in Q compartment die at the rate $\delta \alpha$, where $\delta$ is the rate at which people die. Here, $\Omega$ and $\mu$ are the natural birth rate and natural death rate respectively. We assume that $\Omega = \mu$. 

All the parameters for this work are chosen in the Indian scenario from \cite{main1,main2,main3} making the results applicable in Indian context.

\begin{table}[htp!]
\begin{center}
\begin{tabular}{ | c | c | c | c | c | }
\hline
 \textbf{Variable} & \textbf{Description} & \textbf{Value} & \textbf{Source}\\
 \hline
 $\epsilon$ & Rate of infection & 0.1961 & \cite{main1} \\
 \hline
 $\frac{1}{\gamma}$ & Infectious period & 7 days & \cite{main1} \\  
 \hline
 $\frac{1}{\rho}$ & Quarantine Period & 14 days & \cite{main1} \\
 \hline
 $\delta$ & Death Rate & 1 & \cite{main2} \\
 \hline
 $\alpha$ & Death probability & 0.43$\%$ & \cite{main1} \\
 \hline
 $\mu$ & Natural Death rate & 0.062 & \cite{main3} \\
\hline
\end{tabular}
\caption{Table describing the parameter values}
\label{parameters table}
\end{center}
\end{table}

We incorporated two delays in our model\cite{main2,delay3}. They are :

\begin{itemize}
\item[1.] Time elapsed between infection and identification($\tau$) : Isolated individuals at time t were infected at time t-$\tau$. This delay can be accounted for both delay in testing the asymptomatic individuals and the delay in isolating infected individuals. 
\item[2.] Delay associated with recovery($\kappa$) : Quarantined individuals stay in quarantine for $\kappa$ days after which they are recovered. We consider it to be 14 days.
\end{itemize} 

Among all the parameters in the proposed model, the transmission probability $\beta$ may depend on temperature\cite{main1,temp2,temp3,china}. We consider $\beta$ as a linear or quadratic function of temperature in similar lines with \cite{main1,temp3} as follows : \\
$$\beta(T) = \beta_{0} + \beta_{1} T , \beta_{0} = 0.84, \beta_{1} = -0.00425 $$ 
$$\beta(T) = \beta_{0} - \beta_{1} ( T - T_{M})^2) , \beta_{0} = 0.792, \beta_{1} = 0.000345, T_{M} = 7.73$$

Very few infected individuals at time 0 spreads the virus. So, we take initial value as \\ (S(t),E(t),I(t),Q(t),R(t),D(t)) = (0.999,0,0.001,0,0,0). \\ 

The detection probability ($\sigma$) ranges from 0 to 1 depending on the effectiveness of the system. For the preliminary analysis, we shall consider $\tau$ as 4 days\cite{delay3} and p as 0.4 in the further sections. \\
\section{Mathematical Analysis}
\label{sec3}

\subsection{Reproduction Number($\mathcal{R}_{0}$) : } The basic reproductive number is the mean number of secondary cases that a typical infected case will cause in a population with no immunity to the disease in the absence of interventions to control the infection. 

The variational matrix of the model computed at the infection free state ( S = 1 , E = I = Q = R = D = 0) gives the non-negative new infection F-matrix and non-singular transition V-matrix respectively as follows :
$$ F = \begin{bmatrix} 0 & \beta \\ 0 & 0 \\ \end{bmatrix}    ,  V = \begin{bmatrix} \epsilon + \mu & 0 \\  -\epsilon & \mu + \gamma + p e^{- \gamma \tau} \\ \end{bmatrix} $$

The basic reproduction number $\mathcal{R}_{0}$ = $\rho(FV^{-1})$, where $\rho(FV^{-1})$ is the spectral radius for a next generation matrix $FV^{-1}$ \cite{main} is given as :

$$ \mathcal{R}_{0} = \frac{\beta \epsilon}{(\epsilon + \mu) (\gamma + p e^{-\gamma \tau} + \mu)}$$

\begin{figure}[hbt!]
\centering
\includegraphics[width=5in, height=3in, angle=0]{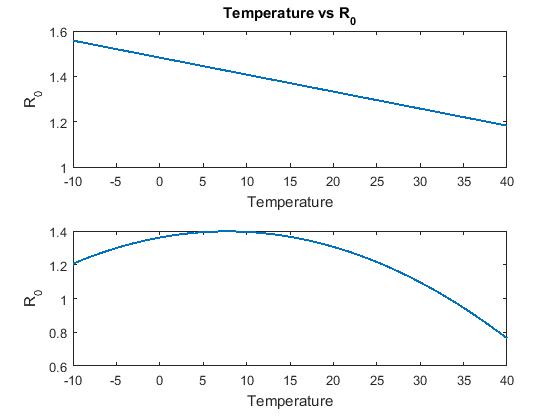}
\caption{Impact of temperature(T) on $\mathcal{R}_{0}$}\label{r}
\end{figure}

\begin{figure}[hbt!]
\centering
\includegraphics[width=5in, height=3in, angle=0]{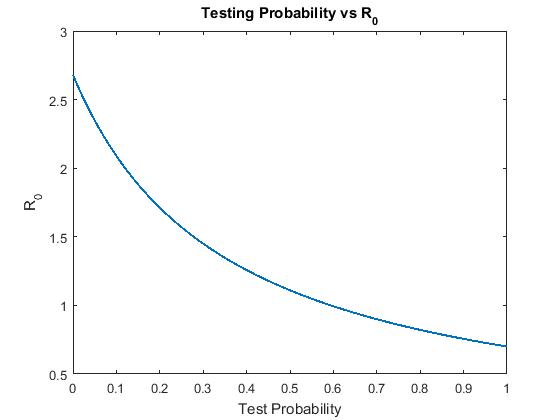}
\caption{Impact of isolation probability(p) on $\mathcal{R}_{0}$}\label{rp}
\end{figure}

\begin{figure}[hbt!]
\centering
\includegraphics[width=5in, height=3in, angle=0]{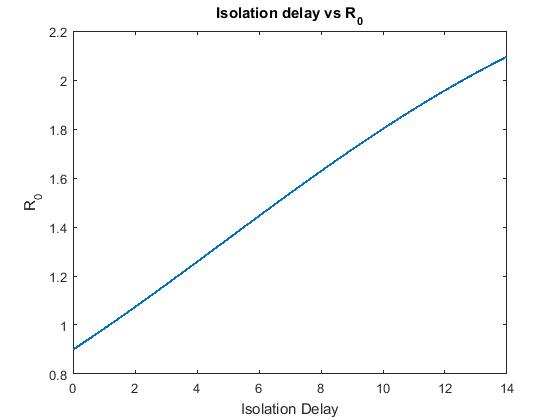}
\caption{Impact of isolation delay($\tau$) on $\mathcal{R}_{0}$}\label{ri}
\end{figure}

The effect of $\beta, p, \tau$ on $\mathcal{R}_{0}$ is discussed in Figures \ref{r},\ref{rp},\ref{ri}. These plots are simulated using dde23 method in MATLAB.

\subsection{Equilibrium states: }

There are two equilibrium points of the delay model which are as follows,

The trivial equilibrium point, E$_{0}$ = (1, 0, 0, 0, 0, 0). \\

The endemic equilibrium is given as :
\begin{align*}
S^{*} &= \frac{1}{\beta \epsilon} (\epsilon + \mu) (\gamma + p e^{-\gamma \tau} + \mu)  \\
E^{*} &= \frac{\mu - \mu S^{*}}{\epsilon + \mu}  \\
I^{*} &= \frac{\epsilon}{\mu + \gamma + p e^{-\gamma \tau}} E^{*} \\
Q^{*} &= \frac{p e^{-\gamma \tau } }{\rho(1-\alpha) + \delta \alpha + \mu} I^{*} \\
R^{*} &= \frac{\gamma I^{*} + \rho(1-\alpha) Q^{*} }{\mu} \\
D^{*} &= \frac{\delta \alpha}{\mu} Q^{*} 
\end{align*} 

\paragraph{Condition for existence of E$^{*}$ }
This equilibrium exists when it satisfies all the initial conditions. So Endemic equilibrium exists if and only if $E^{*} >$ 0. This implies that $S^{*} <$ 1. Hence the condition is $(\epsilon + \mu) (\gamma + p e^{-\gamma \tau} + \mu) < \beta \epsilon$.

\section{Stability Analysis}
\label{sec4}

To perform the stability analysis, the Jacobian matrix[J] for the proposed system is computed as : 
$$J[E] = \begin{bmatrix} -\beta I - \mu & 0 & -\beta S & 0 & 0 & 0 \\ \beta I & -\mu - \epsilon & \beta S & 0 & 0 & 0 \\ 0 & \epsilon & -\gamma - \mu - p e^{-\gamma \tau} e^{-\lambda \tau} & 0 & 0 & 0 \\ 0 & 0 & p e^{-\gamma \tau} e^{-\lambda \tau} & -\rho (1-\alpha)e^{-\kappa \lambda} - \delta \alpha - \mu & 0 & 0 \\ 0 & 0 & \gamma &  -\rho (1-\alpha)e^{-\kappa \lambda} & - \mu & 0 \\ 0 & 0 & 0 & \delta \alpha & 0 & -\mu \end{bmatrix}$$

The characteristic equation is given as :
$$\chi(\lambda) = (\lambda + \mu)^2 (\lambda + \mu + \delta \alpha + \rho (1-\alpha) e^{-\kappa \lambda}) (\lambda^3 + a_1 \lambda^2 + a_2 \lambda + a_3 + e^{-\lambda \tau} (b_0 \lambda^2 + b_1 \lambda + b_2)) $$

where,
\begin{align*}
a_1 &= 3 \mu + \epsilon + \gamma + \beta I \\
a_2 &= (\gamma + \mu) (\epsilon + \mu) + (\gamma + \epsilon + 2 \mu) (\beta I + \mu) - \epsilon \beta S \\
a_3 &= (\gamma + \mu) (\epsilon + \mu) (\beta I + \mu) - \beta \epsilon \mu S \\
b_0 &= p e^{-\gamma \tau} \\
b_1 &= p e^{-\gamma \tau} (2 \mu + \epsilon + \beta I) \\
b_2 &= p e^{-\gamma \tau} (\epsilon + \mu) (\mu + \beta I)
\end{align*}

\subsubsection{Stability of Infection-free state(E$_{0}$) : }

The characteristic equation about the infection-free state(E$_{0}$) is given as :
\begin{equation}
\chi(\lambda) = (\lambda + \mu)^3  (\lambda + \mu + \delta \alpha + \rho (1-\alpha) e^{-\kappa \lambda}) (\lambda^2 + d_1 \lambda + d_2 + e^{-\lambda \tau} ( e_1 \lambda + e_2)) \label{**} \end{equation}

where,
\begin{align*}
d_1 &= 2 \mu + \epsilon + \gamma \\
d_2 &= (\gamma + \mu) (\epsilon + \mu) - \epsilon \beta \\
e_1 &= p e^{-\gamma \tau} \\
e_2 &= p e^{-\gamma \tau} (\mu + \epsilon) 
\end{align*}

Now the characteristic equation without delays $(\tau=0, \kappa=0)$ is given as,
\begin{equation}
\chi(\lambda) = (\lambda + \mu)^3  (\lambda + \mu + \delta \alpha + \rho (1-\alpha)) (\lambda^2 + d_1 \lambda + d_2 +  ( e_1 \lambda + e_2)) \label{1}
\end{equation}
Using Descartes rule of sign change we see that all the  roots of (\ref{1}) are negative, which implies that the disease free equilibrium $E_0$ is locally asymptotically stable with $\alpha < 1$ and $(\gamma+\mu)(\epsilon+\mu) > \epsilon \beta$.\\

Now we will analyze the stability of $E_0$ with delays. Let us first consider the following,

\begin{equation}
    (\lambda + \mu)^3  (\lambda + \mu + \delta \alpha + \rho (1-\alpha) e^{-\kappa \lambda}) \label{*}
\end{equation}
The three roots of this equation are : $\lambda_{1,2,3} = -\mu$, and others are the roots of the following equation.
\begin{equation}
     (\lambda + \mu + \delta \alpha + \rho (1-\alpha) e^{-\kappa \lambda}) \label{2}
\end{equation}

Let $\lambda = i \omega$. Substituting in the above equation (\ref{2}) and simplifying the terms and equating the real parts and imaginary parts, we get, 

\begin{align*}
\omega &= \rho(1-\alpha cos\kappa \omega) \\
\mu+\delta \alpha &= \rho(1-\alpha sin\kappa \omega)
\end{align*}

Squaring and adding the above two equations give
\begin{equation}
\omega^2 + (\mu + \delta \alpha)^2 -\rho^2(1-\alpha)^2 = 0 \label{3}
\end{equation}
If $(\mu + \delta \alpha)^2$ $>$ $\rho^2(1-\alpha)^2$
then all roots of (\ref{3}) are imaginary, which implies that all roots of (\ref{2}) are negative. So far now we have all the roots of (\ref{*}) to be negative provided $(\mu + \delta \alpha)^2  >$ $\rho^2(1-\alpha)^2$.

Now we consider the following expression of (\ref{**}) and state and prove the bifurcation theorem for $E_0$.
\begin{equation}
 ( \lambda^2 + d_1 \lambda + d_2 + e^{-\lambda \tau} (e_1 \lambda + e_2)) = 0 \label{4}
\end{equation}

\paragraph{\textbf{Theorem}}
\label{theorem}
Let $(\mu + \delta \alpha)^2 > (1-\alpha)^2 \rho^2$ and $f_1 > 0$ and $E_0$ be asymptotically stable without any delay. Then \vspace{2mm}
\begin{itemize}
\item[1.] When the conditions $d_2^2 < e_2^2$, $\gamma^2 > \omega^2$, $d_2 < \omega^2$ and $x > y$ holds, the disease free equilibrium point $E_0^*$ of the delayed model (2.1) is locally asymptotically stable while $\tau < \tau^*$ and is unstable while $\tau >\tau^*$ where $$\tau^* = \frac{1}{\omega_0} arccos \bigg[\frac{(e_2 - d_1 e_1) \omega^2 - d_2 e_2}{e_1^2 \omega^2 + e_2^2 }\bigg] $$

For $\tau \;=\; \tau*$, Hopf Bifurcation occurs as the value of $\tau$ passes through the critical value $\tau^*$. \\
\item[2.] If $d_2^2 > e_2^2$, then $E_0$ is locally asymptotically stable for all $\tau.$ \\
\end{itemize}

\noindent\textbf{Proof:}
The disease free equilibrium point $E_0$ of the delayed system (2.1) will be asymptotically stable if all the roots of the characteristic equation (\ref{4}) have negative real parts. We have seen that $E_0$ of the non delayed system is asymptotically stable. We now check the sign of  the real part of the roots of (\ref{4})  as we vary the value of the gestation delay $\tau$.\\

\noindent
Let $\lambda \;=\; \mu(\tau) + i\omega(\tau)$ where, $\mu$ and $\omega$ are real. Since $E_0$ is stable for the non-delayed system and $\mu(0) < 0$, we will choose $\tau$ sufficiently close to 0 and use continuity of $\tau$ to prove the theorem. Let $\tau > 0$ be sufficiently small, then by continuity $ \mu(\tau) < 0 $ and $E^*$ will still remain stable. The stability changes for some values of $\tau$ for which $\mu(\tau) \;=\; 0 $ and $\omega(\tau)\ne 0 $ that is when $\lambda$ is purely imaginary. Let $\tau^*$ be such that $\mu(\tau^*) \;=\;0 $ and $\omega(\tau^*)\ne 0$. In this case the steady state loses stability and becomes unstable when $\mu(\tau)$ becomes positive.\\

By Rouche's theorem \citep*{hi} the transcendental equation (\ref{**}) has roots with positive real parts if and only if it has purely imaginary roots. Now we will assume that the characteristic equation (\ref{4}) has purely imaginary roots and then arrive at a contradiction.\\  
Let $\lambda = i \omega$. Substituting in the above equation, simplifying the terms and equating the real parts and imaginary parts, we get, 

\begin{align*}
d_2 - \omega^2 &= (- e_1 \omega) sin(\omega \tau) - e_2 cos(\omega \tau) \\
d_1 \omega &= (-e_1 \omega) cos(\omega \tau) + e_2 sin(\omega \tau)
\end{align*}

Squaring and adding the above two equations give
\begin{equation}
\omega^4 + f_1 \omega^2 + f_2 = 0 \label{omega}
\end{equation}
where, 
\begin{align*}
f_1 &= d_1^2 - 2 d_2 - e_1^2 \\
f_2 &= d_2^2 - e_2^2 
\end{align*}

If $d_2^2 > e_2^2$ then all roots of (\ref{omega}) are imaginary. This contradicts that $\omega$ is real therefore equation (\ref{4}) does not have purely imaginary roots. Hence by Rouche's Theorem we conclude that all roots of (\ref{4}) have negative real parts. Therefore the disease free equilibrium point $E_0$ is asymptotically stable for all $\tau$, whenever $d_2^2 > e_2^2 $. 

When $d_2^2 < e_2^2$ from equation (\ref{omega}), we see that there is only one sign change, therefore using Descartes rule of sign  change we conclude that whenever $d_2^2 < e_2^2$  equation (\ref{omega}) has unique positive root say $\omega^*$. This implies that the characteristics equation (\ref{4}) has a purely imaginary roots $\pm i\omega^*$. Now we will calculate the critical value of the delay $\tau^*$ based on $\omega^*.$ \vspace{5mm}\\

\noindent
\underline{Calculation of $\tau^*$ at $\omega^*$:} \\
Consider,
\begin{align*}
d_2 - \omega^{*2} &= (- e_1 \omega^*) sin(\omega^* \tau) - e_2 cos(\omega^* \tau) \\
d_1 \omega^* &= (-e_1 \omega^*) cos(\omega^* \tau) + e_2 sin(\omega^* \tau)
\end{align*}
Solving the above equation for $\tau^*$ we get,
\begin{equation}
\tau^* = \frac{1}{\omega^*} cos^{-1} \bigg[\frac{(e_2 - d_1 e_1) \omega^2 - d_2 e_2}{e_1^2 \omega^2 + e_2^2 }\bigg] 
\end{equation}

Now to show that there is at least one eigen value with a positive real part for $\tau > \tau^*$, we verify the transversality condition $\frac{d(Re \lambda(\tau))}{d\tau} > 0 $ at $\tau\;=\;\tau^*$ \\

From (\ref{4}) we have, \\
$$ \lambda^2 + d_1 \lambda + d_2 + e^{-\lambda \tau}(e_1 \lambda + e_2)=0$$
substituting for $e_1$ and $e_2$ we get
\begin{equation}
 \lambda^2 + d_1 \lambda + d_2 + (\lambda \rho + \rho(\mu + \epsilon))e^{-(\lambda + \gamma)\tau} =0 \label{c}
 \end{equation}
\noindent
Differentiating the above equation with respect to $\tau,$ we get, \\

$$\bigg(2\lambda+ d_1 + \rho e^{-(\lambda + \gamma)\tau}- \lambda \rho \tau(\lambda +  \gamma) e^{-(\lambda+\gamma)\tau}-\rho \tau (\mu + \epsilon)(\lambda +\gamma) e^{-(\lambda+\gamma)\tau}\bigg)\frac{d \lambda}{d \tau}$$ $$= \lambda \rho (\lambda + \gamma)^2 e^{-(\lambda + \gamma)\tau} + \rho (\mu + \epsilon)(\lambda + \gamma)^2 e^{-(\lambda + \gamma)}$$

simplifying we get, 
\begin{equation}
\hspace{1cm}\bigg(\frac{d\lambda}{d\tau}\bigg)^{-1} \;=\;\frac{2\lambda+d_1}{\rho (\lambda + \gamma)^2(\lambda + \mu + \epsilon)e^{-(\lambda + \gamma)\tau}} + \frac{1-\tau(\lambda + \gamma)(\lambda + \mu + \epsilon)}{ (\lambda + \gamma)^2(\lambda + \mu + \epsilon)}\label{b}
\end{equation}
 From (\ref{c}) we get,
 $$e^{-(\lambda +\gamma)\tau} = \frac{-\lambda^2-d_1 \lambda-d_2}{\rho(\mu+\epsilon)+\lambda \rho}$$
  
  substituting the expression for $e^{-(\lambda +\gamma)\tau}$  in (\ref{b}) we get,
  $$\bigg(\frac{d\lambda}{d\tau}\bigg)^{-1}=\frac{-(2\lambda+d_1)}{(\lambda+\gamma)^2(\lambda^2+d_1 \lambda + d_2)}+ \frac{1}{(\lambda+\gamma)^2(\lambda+\mu+\epsilon)}-\frac{\tau}{(\lambda+\gamma)}$$
  
At $\tau\;=\;\tau^*$ and $\lambda \;=\;i\omega^*$
$$\bigg(\frac{d\lambda}{d\tau}\bigg)^{-1}=\frac{-(2i\omega+d_1)}{(i\omega+\gamma)^2((i\omega)^2+d_1 i\omega + d_2)}+ \frac{1}{(i\omega+\gamma)^2(i\omega+\mu+\epsilon)}-\frac{\tau}{(i\omega+\gamma)}$$
Now considering only the real parts and simplifying we get,

$$\bigg(Re(\frac{d\lambda}{d\tau})^{-1}\bigg)=\frac{d_1(-\gamma^2+\omega^2)(d_2-\omega^2)+ 2 \omega^2 d_1 \gamma(2+d_1)}{((\gamma^2-\omega^2)^2+4 \omega^2 \gamma^2)((d_2-\omega^2)^2+d_1\omega^2)}+\frac{(\gamma^2-\omega^2)(\mu+\epsilon)-2\omega^2\gamma}{((\gamma^2-\omega^2)^2+4\omega^2\gamma^2)((\mu+\epsilon)^2+\omega^2)}-\frac{\tau^* \gamma}{(\gamma^2+\omega^2)}$$

Therefore, we have,

$$sign\bigg(Re(\frac{d\lambda}{d\tau})^{-1}\bigg) > 0 $$ 
provided we have $(x + y )> 0$ and $x + y > z $, where,
$$x=\frac{d_1(-\gamma^2+\omega^2)(d_2-\omega^2)+ 2 \omega^2 d_1 \gamma(2+d_1)}{((\gamma^2-\omega^2)^2+4 \omega^2 \gamma^2)((d_2-\omega^2)^2+d_1\omega^2)}$$
$$y=\frac{(\gamma^2-\omega^2)(\mu+\epsilon)-2\omega^2\gamma}{((\gamma^2-\omega^2)^2+4\omega^2\gamma^2)((\mu+\epsilon)^2+\omega^2)}$$
$$z=\frac{\tau^* \gamma}{(\gamma^2+\omega^2)}$$
 Hence with the above conditions transversality Condition holds. Therefore this implies that there is at least one eigenvalue with the positive real part for $\tau > \tau^*$. Hence Hopf Bifurcation occurs at $\tau=\tau^*$. 
\subsubsection{\textbf{Stability of Endemic state(E$^{*}$) : }} 

The characteristic equation about the infection-free state(E$^{*}$) is given as :
$$\chi(\lambda) = (\lambda + \mu)^2 (\lambda + \mu + \delta \alpha + \rho (1-\alpha) e^{-\kappa \lambda}) (\lambda^3 + a_1 \lambda^2 + a_2 \lambda + a_3 + e^{-\lambda \tau} (b_0 \lambda^2 + b_1 \lambda + b_2)) $$

where,
\begin{align*}
a_1 &= 3 \mu + \epsilon + \gamma + \mu (\mathcal{R}_{0} -1) \\
a_2 &= (\gamma + \epsilon + 2 \mu) \mu \mathcal{R}_{0} - p e^{-\gamma \tau} (\epsilon + \mu) \\
a_3 &= \mu (\epsilon + \mu) [ (\gamma + \mu ) ( \mathcal{R}_{0} - 1) - p e^{-\gamma \tau}] \\
b_0 &= p e^{-\gamma \tau} \\
b_1 &= p e^{-\gamma \tau} ( \mu + \epsilon + \mu \mathcal{R}_{0} ) \\
b_2 &= p e^{-\gamma \tau} (\epsilon + \mu) \mu \mathcal{R}_{0}
\end{align*}

Now the characteristic equation without delays $(\tau=0, \kappa=0)$ is given as,
\begin{equation}
\chi(\lambda) = (\lambda + \mu)^2  (\lambda + \mu + \delta \alpha + \rho (1-\alpha)) (\lambda^3 + (a_1+b_0) \lambda^2 + (a_2+b_1)\lambda + a_3+b_2 ( e_1 \lambda + e_2)) \label{cc}
\end{equation}

Using Descartes rule of sign change we see that all the  roots of (\ref{cc}) are negative provided we have,\\
1. $\alpha < 1$\\
2. $(a_1+b_0) > 0$, $(b_1+a_2) > 0$, and $(a_3+b_2) > 0$

Therefore with the above conditions the infected equilibrium $E_1$ to be locally asymptotically stable without any delay.\\

Now we will analyze the stability of $E_0$ with delays. Let us first consider the following,

\begin{equation}
    (\lambda + \mu)^2  (\lambda + \mu + \delta \alpha + \rho (1-\alpha) e^{-\kappa \lambda}) \label{***}
\end{equation}
The three roots of this equation are : $\lambda_{1,2} = -\mu$, and others are the roots of the following equation.
\begin{equation}
     (\lambda + \mu + \delta \alpha + \rho (1-\alpha) e^{-\kappa \lambda}) \label{222}
\end{equation}

Let $\lambda = i \omega$. Substituting in the above equation (\ref{222}) and simplifying the terms and equating the real parts and imaginary parts, we get, 

\begin{align*}
\omega &= \rho(1-\alpha cos\kappa \omega) \\
\mu+\delta \alpha &= \rho(1-\alpha sin\kappa \omega)
\end{align*}

Squaring and adding the above two equations give
\begin{equation}
\omega^2 + (\mu + \delta \alpha)^2 -\rho^2(1-\alpha)^2 = 0 \label{333}
\end{equation}
If $(\mu + \delta \alpha)^2$ $>$ $\rho^2(1-\alpha)^2$
then all roots of (\ref{333}) are imaginary, which implies that all roots of (\ref{222}) are negative. So far now we have all the roots of (\ref{***}) to be negative provided $(\mu + \delta \alpha)^2  >$ $\rho^2(1-\alpha)^2$.

Now consider the following,
\begin{equation}
 (\lambda^3 + a_1 \lambda^2 + a_2 \lambda + a_3 + e^{-\lambda \tau} (b_0 \lambda^2 + b_1 \lambda + b_2)) = 0 \label{sa}
\end{equation}

Let $\lambda = i \omega$. Substituting in the above equation gives :

$$- i \omega^3 - a_1 \omega^2 + a_2 i \omega + a_3 (cos \omega \tau - i sin \omega \tau) (-b_0 \omega^2 + b_1 i \omega + b_2 = 0$$

Simplifying the terms and equating the real parts and imaginary parts, we get, 

\begin{align*}
a_3 - a_1 \omega^2 &= (-a_3 b_1 \omega) sin(\omega \tau) + ( a_3 b_0 \omega^2 - a_3 b_2) cos(\omega \tau) \\
a_2 \omega - \omega^3 &= (-a_3 b_1 \omega) cos(\omega \tau) - ( a_3 b_0 \omega^2 - a_3 b_2) sin(\omega \tau)
\end{align*}

Squaring and adding the above two equations give 

\begin{equation}
\omega^6 + c_1 \omega^4 + c_2 \omega^2 + c_3 = 0 \label{bbc}
\end{equation}
where, 
\begin{align*}
c_1 &= a_1^2 - 2 a_2 - a_3^2 b_0^2 \\
c_2 &= a_2^2 + 2 a_3^2 b_0 b_2 - a_3^2 b_1^2 - 2 a_1 a_3\\
c_3 &= a_3^2 (1 - b_2^2) 
\end{align*}

\begin{figure}[hbt!]
\centering
\includegraphics[width=5in, height=3in, angle=0]{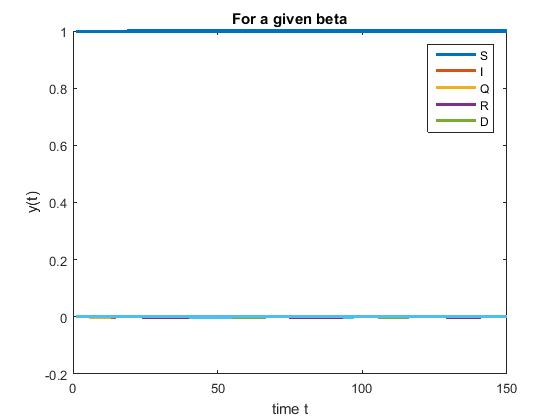}
\caption{ System dynamics when R$_{0} < $1 ($\beta = $ 0.5)}\label{s1}
\end{figure}

\begin{figure}[hbt!]
\centering
\includegraphics[width=5in, height=3in, angle=0]{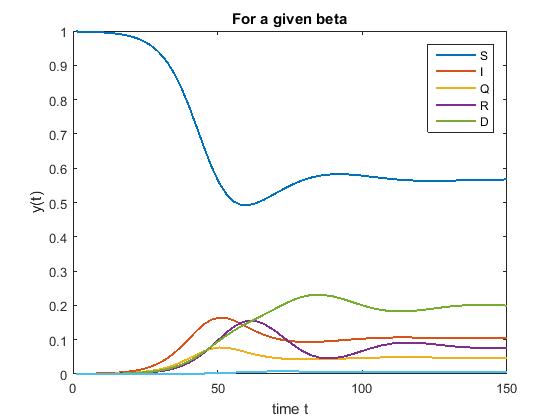}
\caption{System dynamics when R$_{0} > $1 ($\beta = $ 1.0)}\label{s2}
\end{figure}

Now if $c_1 > 0,c_2 > 0, c_3 > 0$ then (\ref{bbc}) has no real roots. This implies that the characteristic equation (\ref{sa}) has all roots to be negative. Therefore $E^*$ is asymptotically stable. Where if $c_3 < 0$ the we have $E^*$ to be unstable.
Thus from the above  analysis we conclude that $E^*$ is asymptotically stable provided:\\
1.$(\mu + \delta \alpha)^2$ $>$ $\rho^2(1-\alpha)^2$\\
2. $c_1 > 0,c_2 > 0, c_3 > 0$ \\
If any of the above two condition is violated then $E^*$ should change its nature of being stable. Figure \ref{s1} and Figure \ref{s2} depicts the dynamics of the system for two specific values of $\beta$. We will illustrate this by performing numerical simulations in section \ref{sec5}. 

Figure \ref{s1} and Figure \ref{s2} depicts the dynamics of the system for two specific values of $\beta$. i.e., 0.5 and 1.

\newpage




\section{Sensitivity Analysis}
\label{sec8}

In this section we check the sensitivity of all the parameters of the proposed model. As each parameter is varied in different intervals, the total infected cell population, mean infected cell population and the mean square error are plotted with respect to time. These plots are used to determine the sensitivity of the parameter. The different intervals chosen are given in the following Table \ref{t}. The fixed parameter values are taken from Table \ref{parameters table}. 

\begin{table}[hbt!]
		\caption{Sensitivity Analysis}
		\centering
		\label{t}
		{
			\begin{tabular}{|l|l|l|}
				\hline
				\textbf{Parameter} & \textbf{Interval} & \textbf{Step Size}  \\
				\hline

				$\mu$ & 0 to 0.5 & 0.01 
					\\ \cline{2-2}
				 & .5 to 2.5  &\\
				 \hline
				 $\beta$ & 0 to .5 & 0.01 
					\\ \cline{2-2}
				 & 2  to 3 &\\
				 
				 \hline
				 	$ \alpha $ & 0 to 2 & 0.01 
					\\ \cline{2-2}
				 & 2 to 5  &  \\
				 \hline
				$\gamma$ & 0 to 1 & 0.01 
				\\ \cline{2-2} 
				& 1 to 2.5 & 
				\\ \hline
				
				$\epsilon$ & 0 to 0.5& 0.01 
				\\ \cline{2-2}
				& 1.5 to 2.5 &  \\ 
				\hline 
					$\omega $ & 0 to 2 & 0.01 
				\\ \cline{2-2}
				& 2 to 4 &  \\ 
				\hline 
				
					$\delta$ & 0 to 1 & 0.01 
				\\ \cline{2-2}
				& 1 to 2.5 &  \\ 					
				\hline 
			\end{tabular}
		}
	\end{table}
	\newpage
	
	\subsubsection{Parameter $\boldsymbol{\beta}$}
	
\begin{figure}[hbt!]
			\begin{center}
				\includegraphics[width=2in, height=1.8in, angle=0]{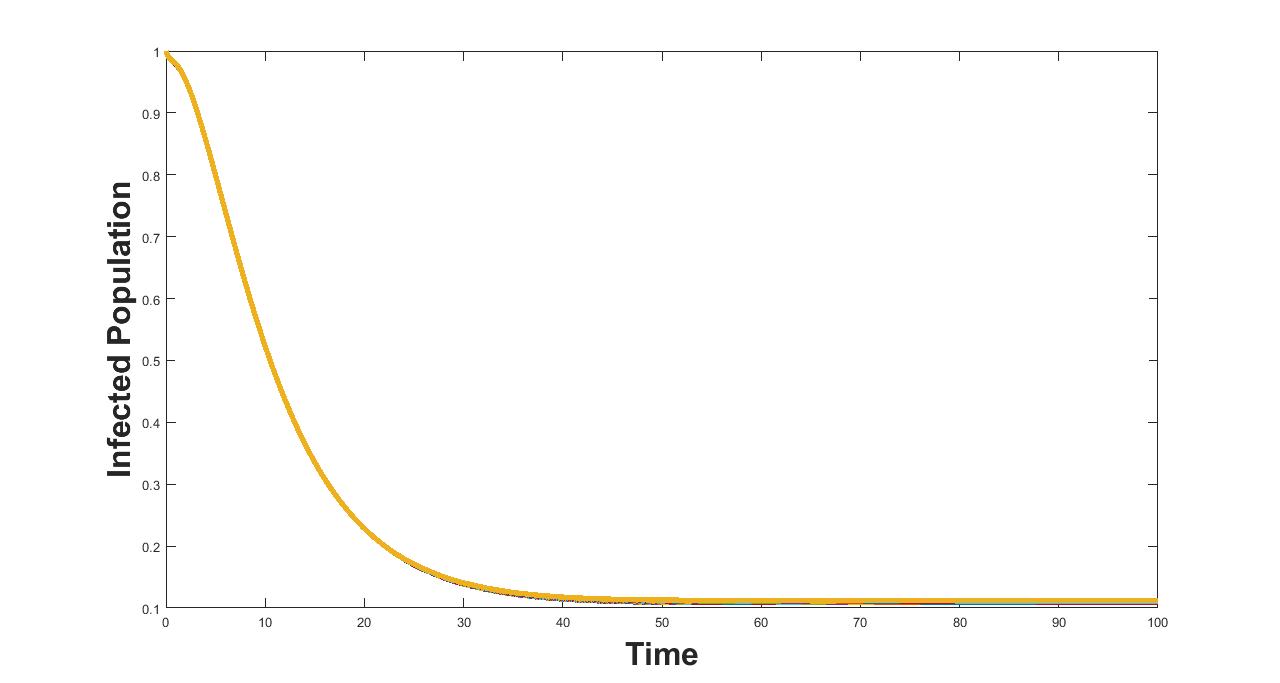}
				\hspace{-.4cm}
				\includegraphics[width=2in, height=1.8in, angle=0]{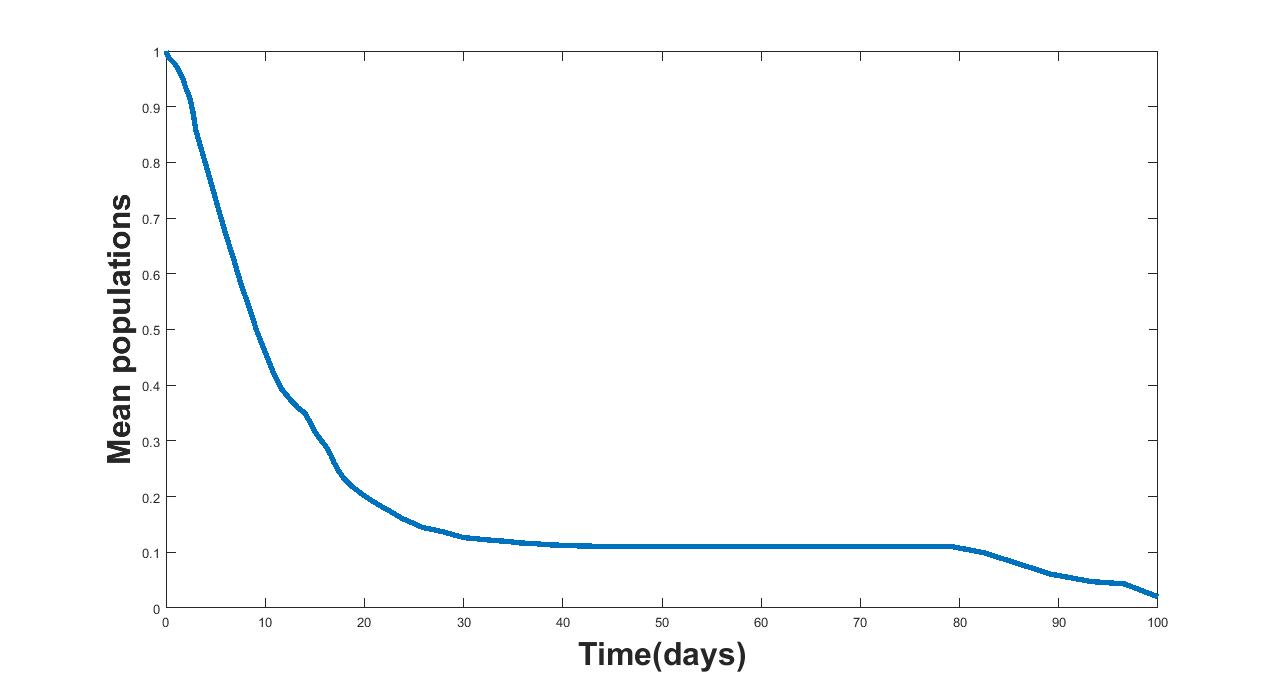}
				\hspace{-.395cm}
					\includegraphics[width=2in, height=1.8in, angle=0]{betamse.jpg}
			\caption*{(b) Interval I :  2 to 3}
				
	 \vspace{7mm}
			\caption{Figure depicting the sensitivity Analysis of $b_1$ varied in three intervals in table \ref{t}. The plots depict the infected population for each varied value of the parameter $\beta$ per interval along with  the mean infected population and the mean square error in the same interval.}
				\label{b1}
			\end{center}
		\end{figure}

	The results related to sensitivity of $\beta$, varied in three intervals as mentioned in table \ref{t}, are given in Figure \ref{b1}. The plots of infected population for each varied value of the parameter $\beta$ per interval, the mean infected population and the mean square error are used to determine the sensitivity.	We conclude from these plots that the parameter $\beta$ is not sensitive in interval II and sensitive in I . 
			\begin{figure}[hbt!]
			\begin{center}
				\includegraphics[width=2in, height=1.8in, angle=0]{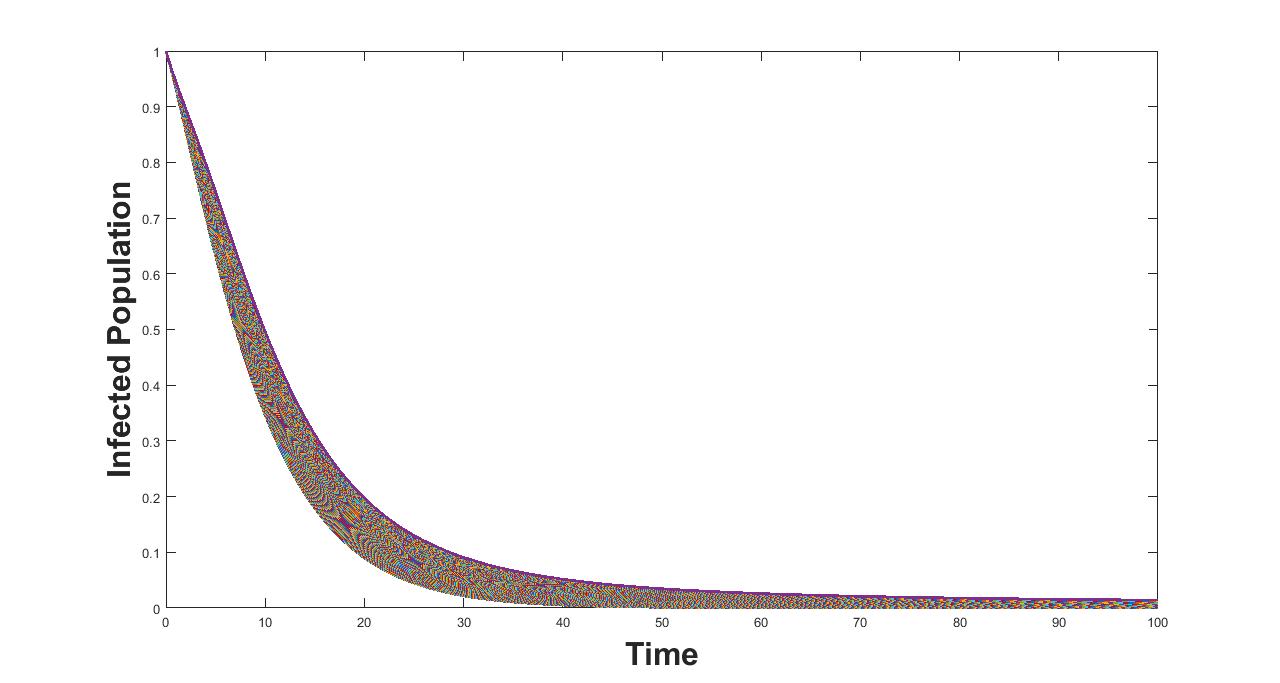}
				\hspace{-.4cm}
				\includegraphics[width=2in, height=1.8in, angle=0]{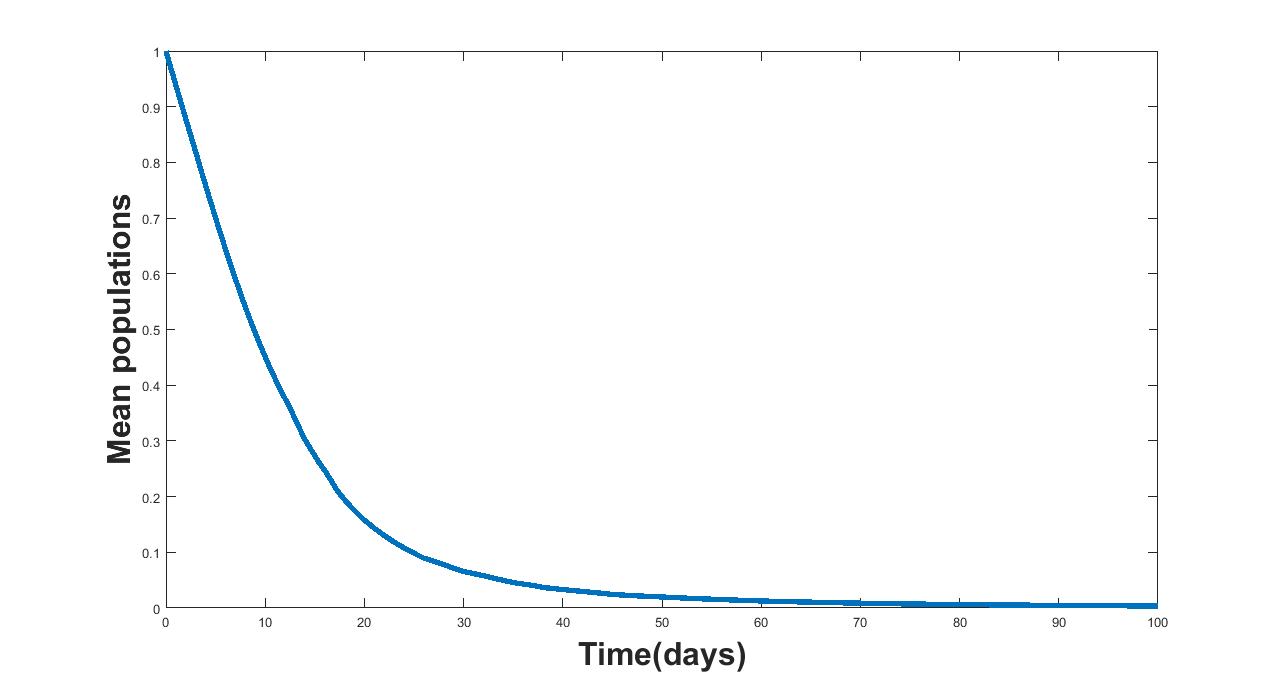}
				\hspace{-.4cm}
					\includegraphics[width=2in, height=1.8in, angle=0]{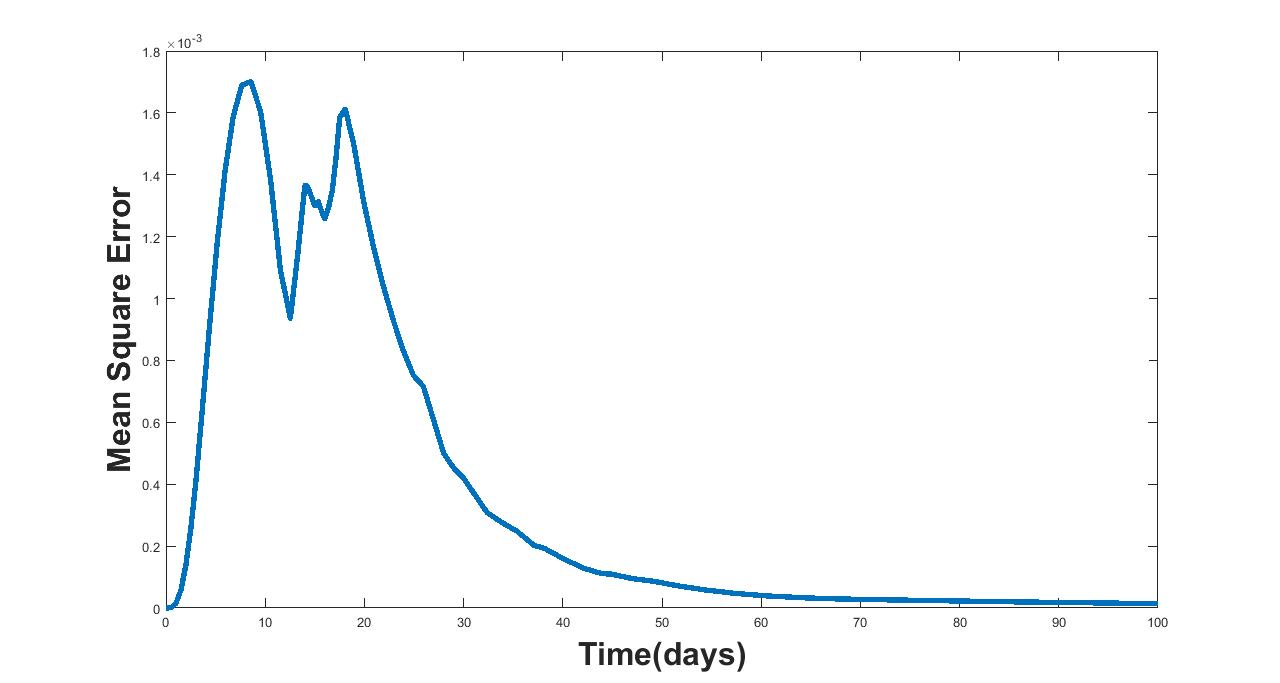}
			\caption*{(a) Interval I :  0 to 0.5}
			\end{center}
		\end{figure}
\newpage
\subsubsection{Parameter $\boldsymbol{\omega}$}

\begin{figure}[hbt!]
			\begin{center}
				\includegraphics[width=2in, height=1.8in, angle=0]{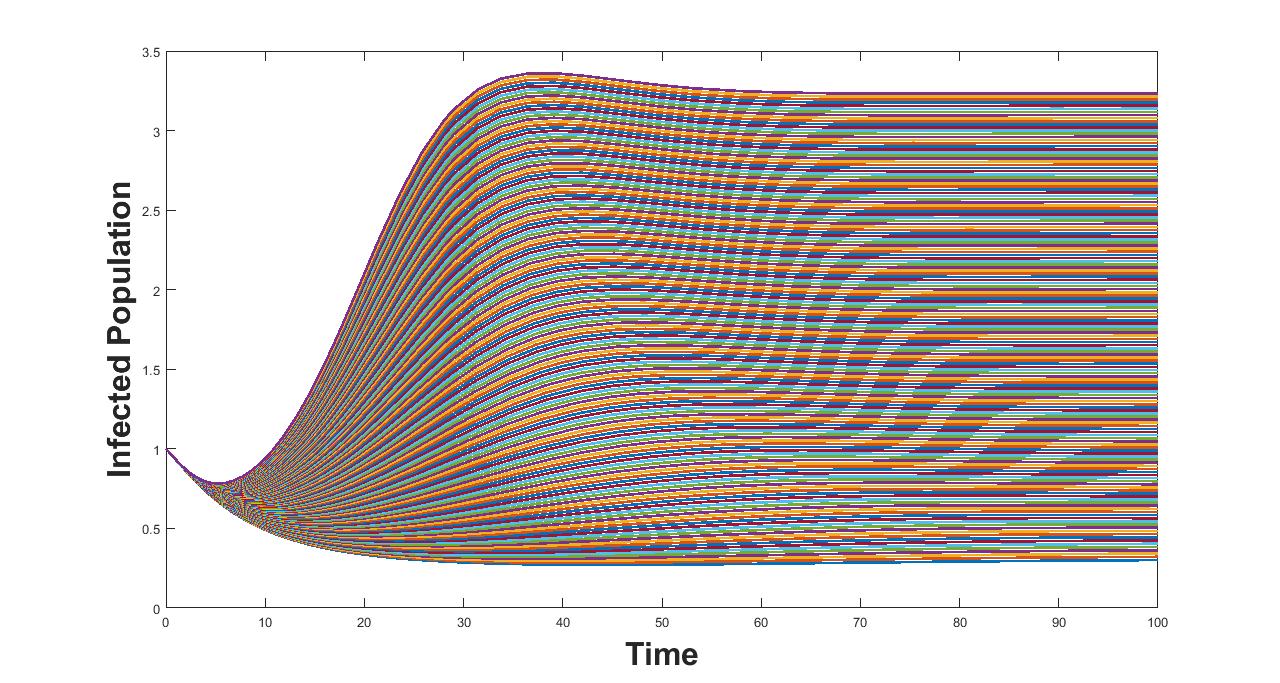}
				\hspace{-.4cm}
				\includegraphics[width=2in, height=1.8in, angle=0]{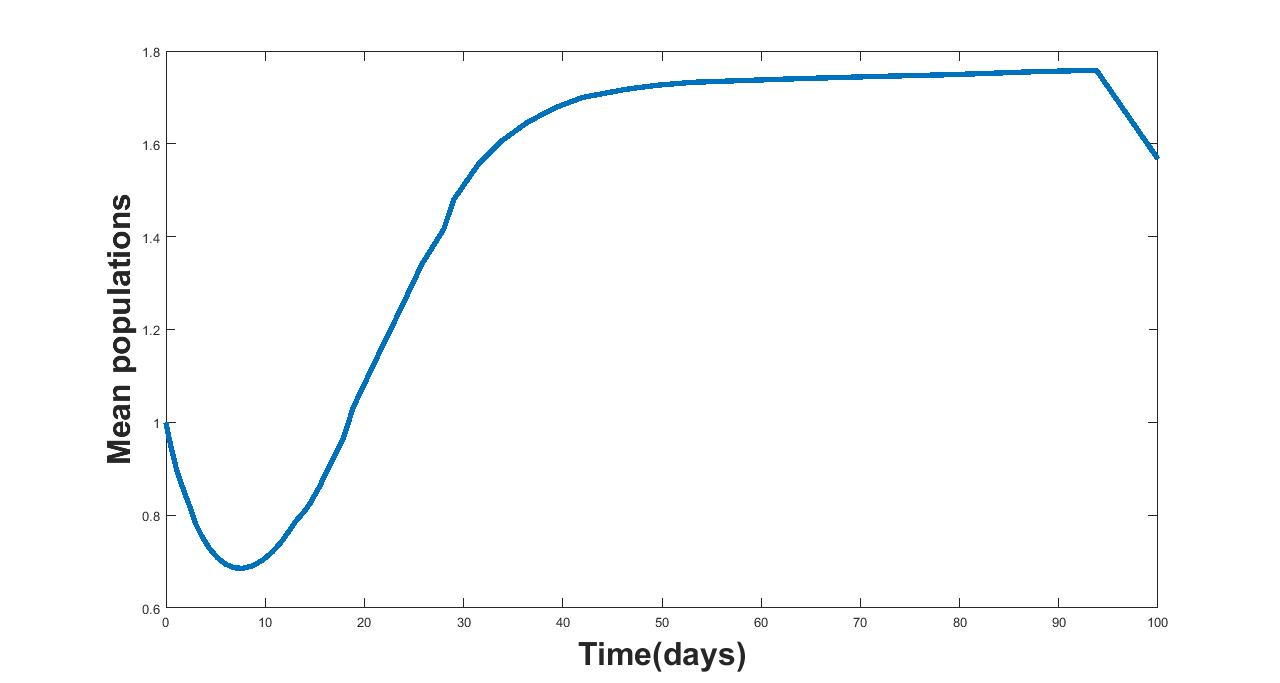}
				\hspace{-.395cm}
					\includegraphics[width=2in, height=1.8in, angle=0]{omegaerrorsen.jpg}
			\caption*{(b) Interval I :  2 to 4}
		
			\vspace{7mm}
			\caption{Figure depicting the sensitivity Analysis of $\omega$ varied in three intervals in table \ref{t}. The plots depict the infected population for each varied value of the parameter $\omega$ per interval along with  the mean infected population and the mean square error in the same interval.    }
				\label{mu}
			\end{center}
		\end{figure}

	The results related to sensitivity of $\omega$, varied in two intervals as mentioned in table \ref{t}, are given in Figure \ref{mu}. The plots of infected population for each varied value of the parameter $\omega$ per interval, the mean infected population and the mean square error are used to determine the sensitivity. We conclude from these plots that the parameter $\omega$ is sensitive in both interval I and II.
		\begin{figure}[hbt!]
			\begin{center}
				\includegraphics[width=2in, height=1.8in, angle=0]{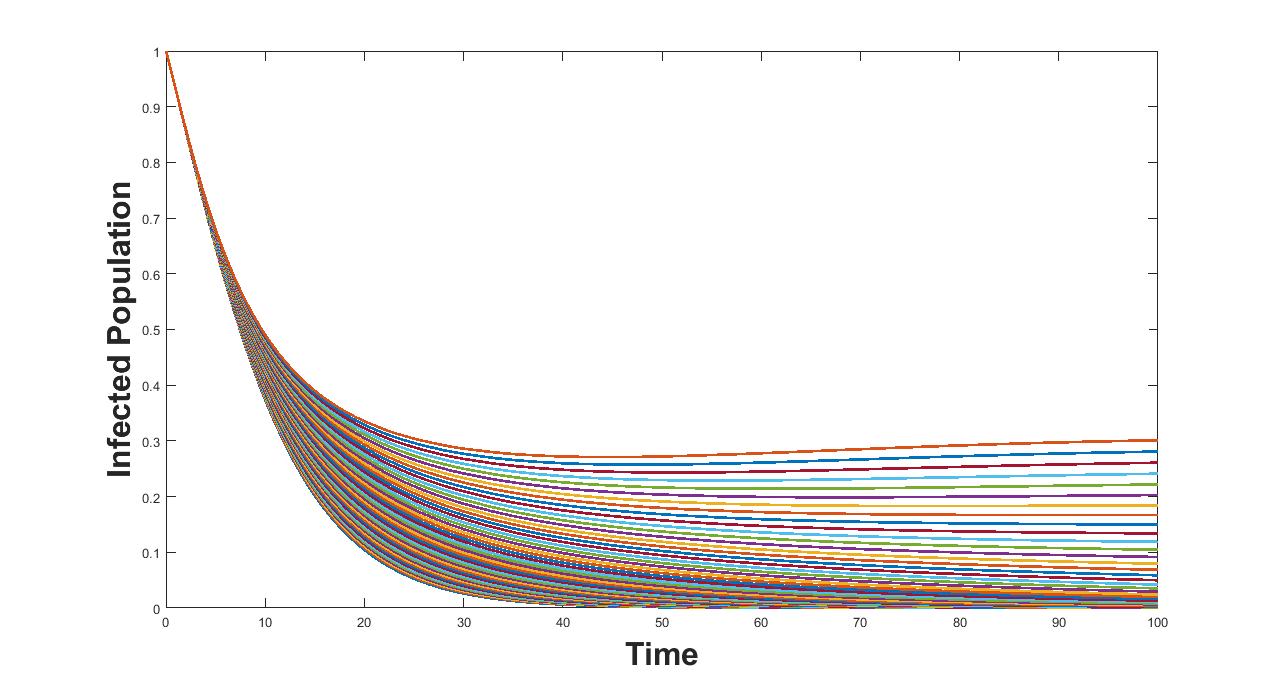}
				\hspace{-.4cm}
				\includegraphics[width=2in, height=1.8in, angle=0]{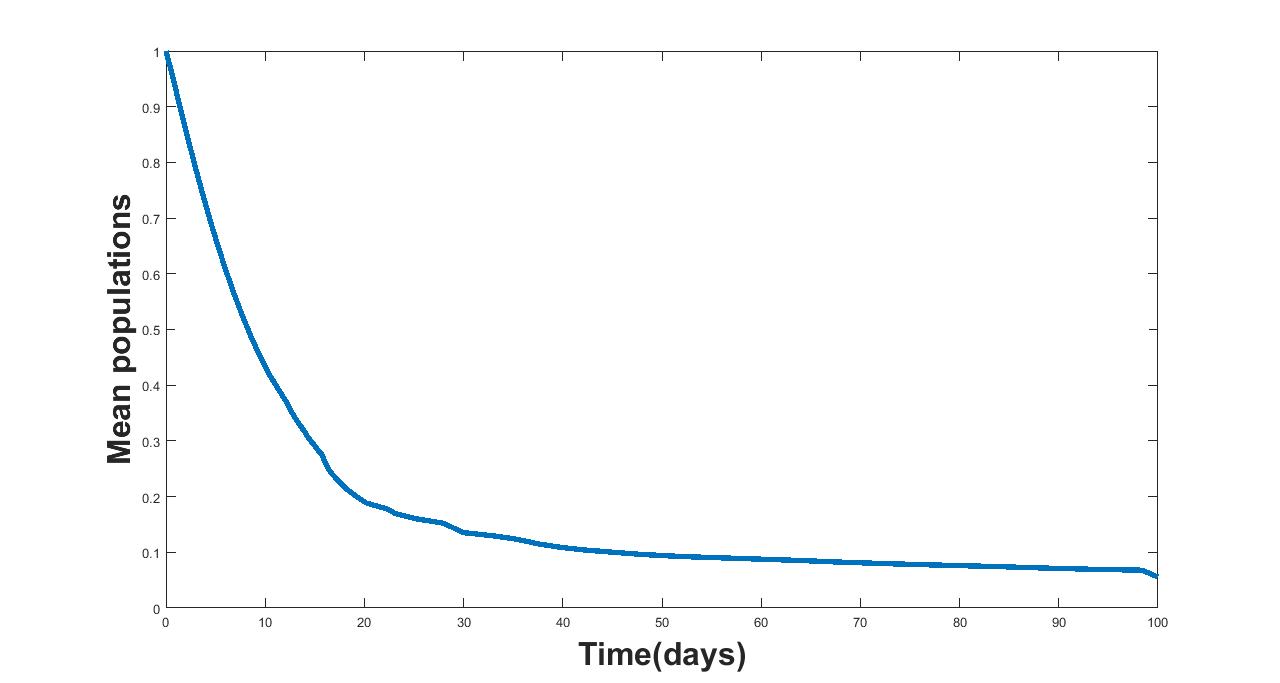}
				\hspace{-.395cm}
					\includegraphics[width=2in, height=1.8in, angle=0]{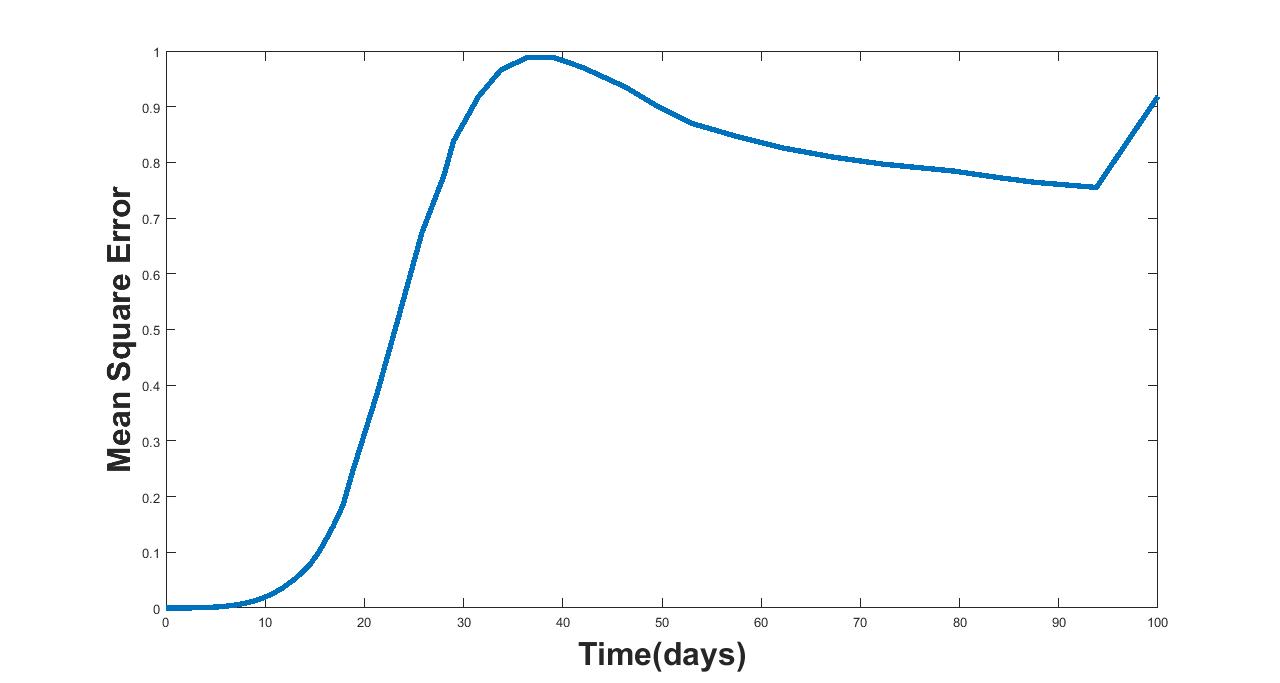}
			\caption*{(a) Interval I :  0 to 2}
			\end{center}
		\end{figure}

		\newpage
		\subsubsection{Parameter $\boldsymbol{\mu}$}

\begin{figure}[hbt!]
			\begin{center}
				\includegraphics[width=2in, height=1.8in, angle=0]{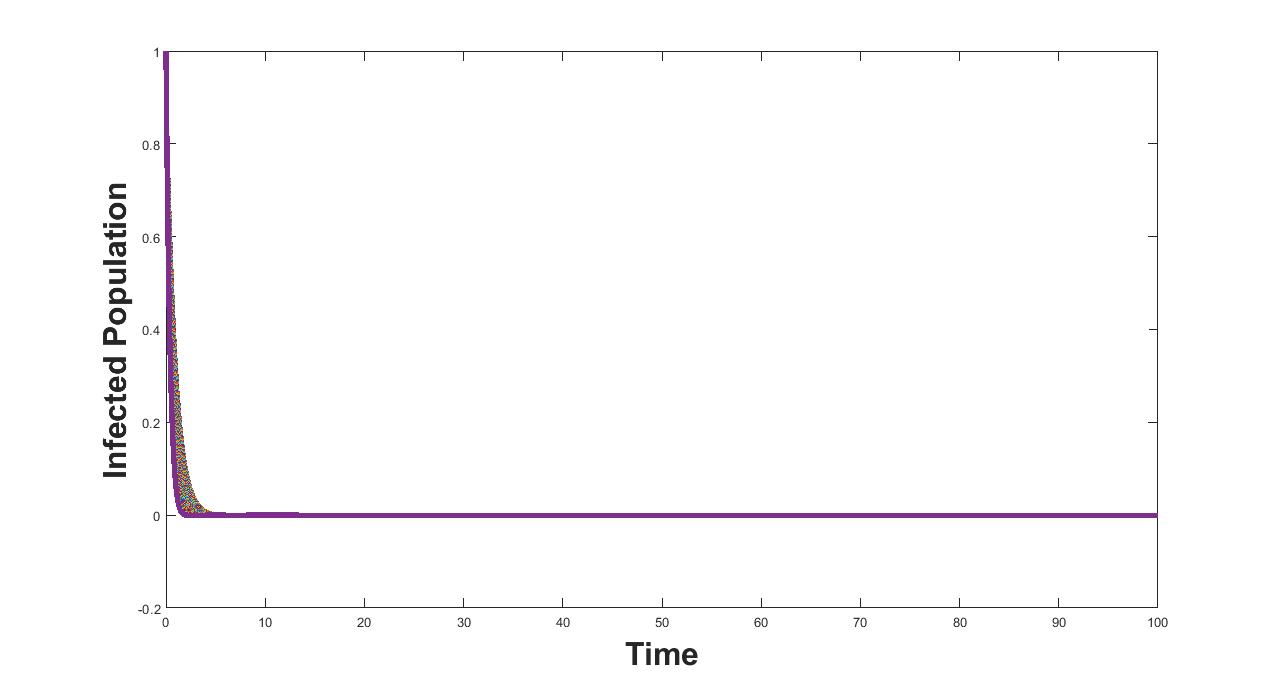}
				\hspace{-.4cm}
				\includegraphics[width=2in, height=1.8in, angle=0]{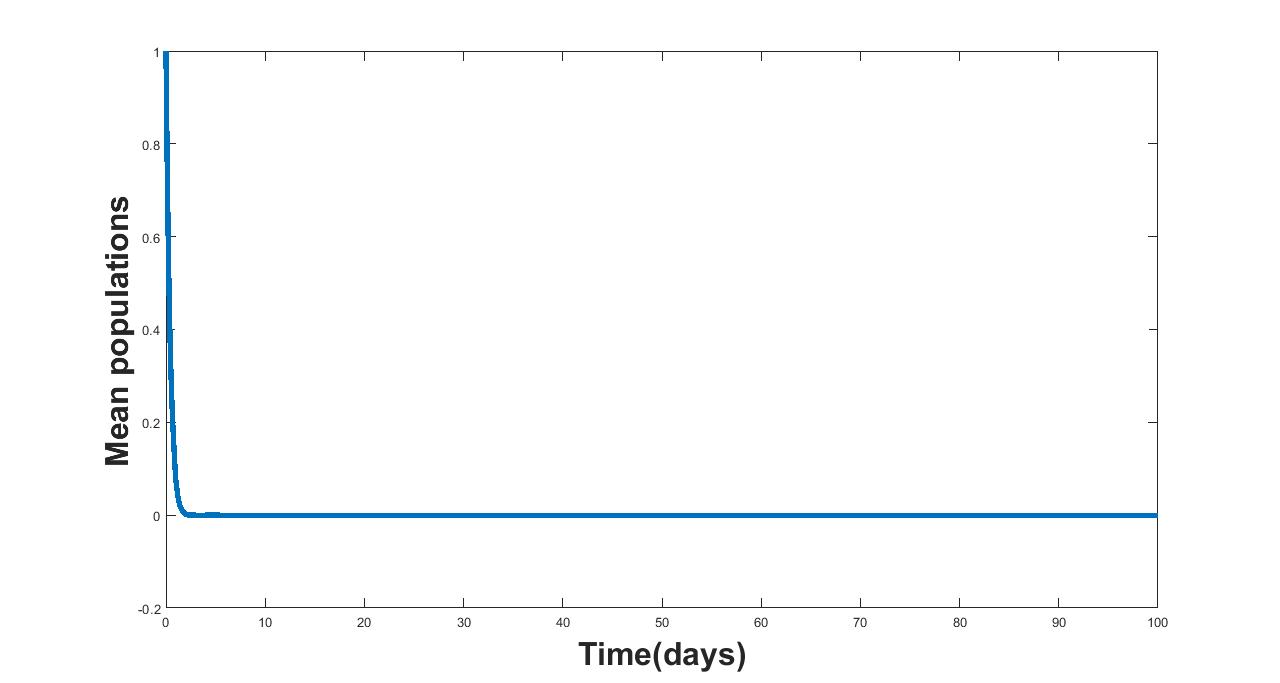}
				\hspace{-.395cm}
					\includegraphics[width=2in, height=1.8in, angle=0]{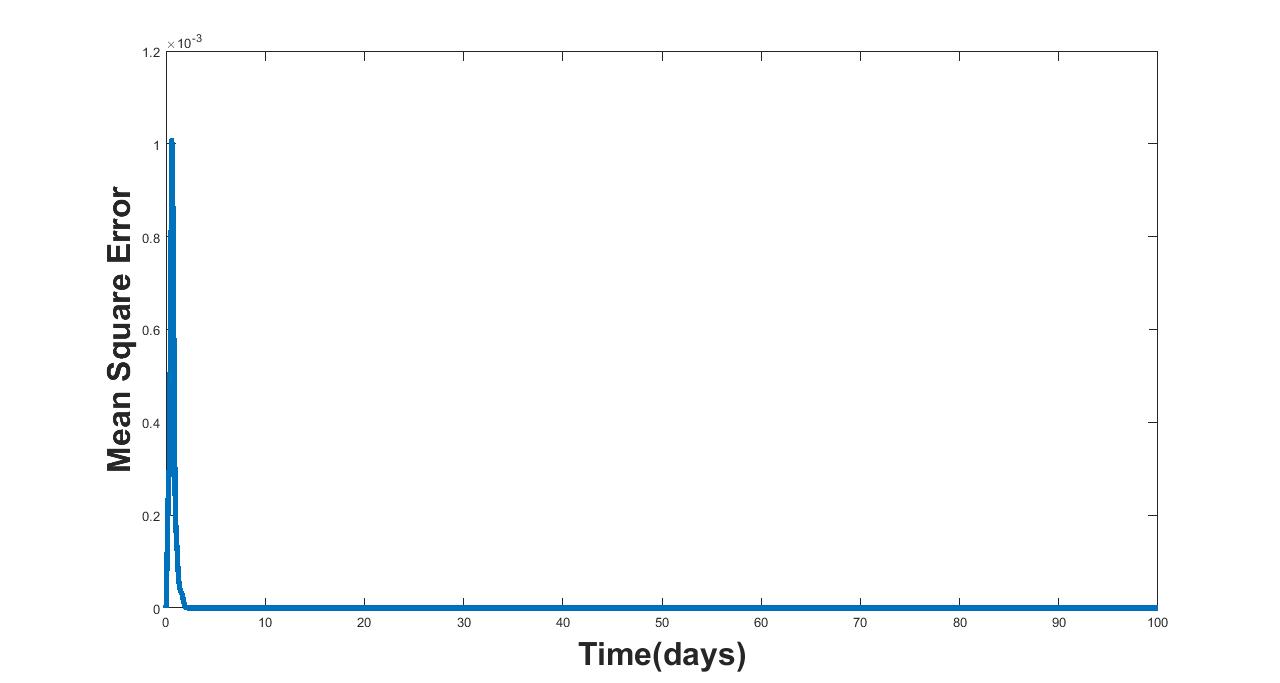}
			\caption*{(b) Interval I :  0.5 to 2.5}
		
			\vspace{7mm}
			\caption{Figure depicting the sensitivity Analysis of $\mu$ varied in two intervals in table \ref{t}. The plots depict the infected population for each varied value of the parameter $\mu$ per interval along with  the mean infected population and the mean square error in the same interval.    }
				\label{beta1}
			\end{center}
		\end{figure}

		The results related to sensitivity of $\mu$, varied in two intervals as mentioned in table \ref{t}, are given in Figure \ref{beta1}. The plots of infected population for each varied value of the parameter $\mu$ per interval, the mean infected population and the mean square error are used to determine the sensitivity. We conclude from these plots that the parameter $\mu$ is sensitive in interval I and insensitive in II. In similar lines, the sensitivity analysis is done for other parameters. 

		\begin{figure}[hbt!]
			\begin{center}
			\includegraphics[width=2in, height=1.8in, angle=0]{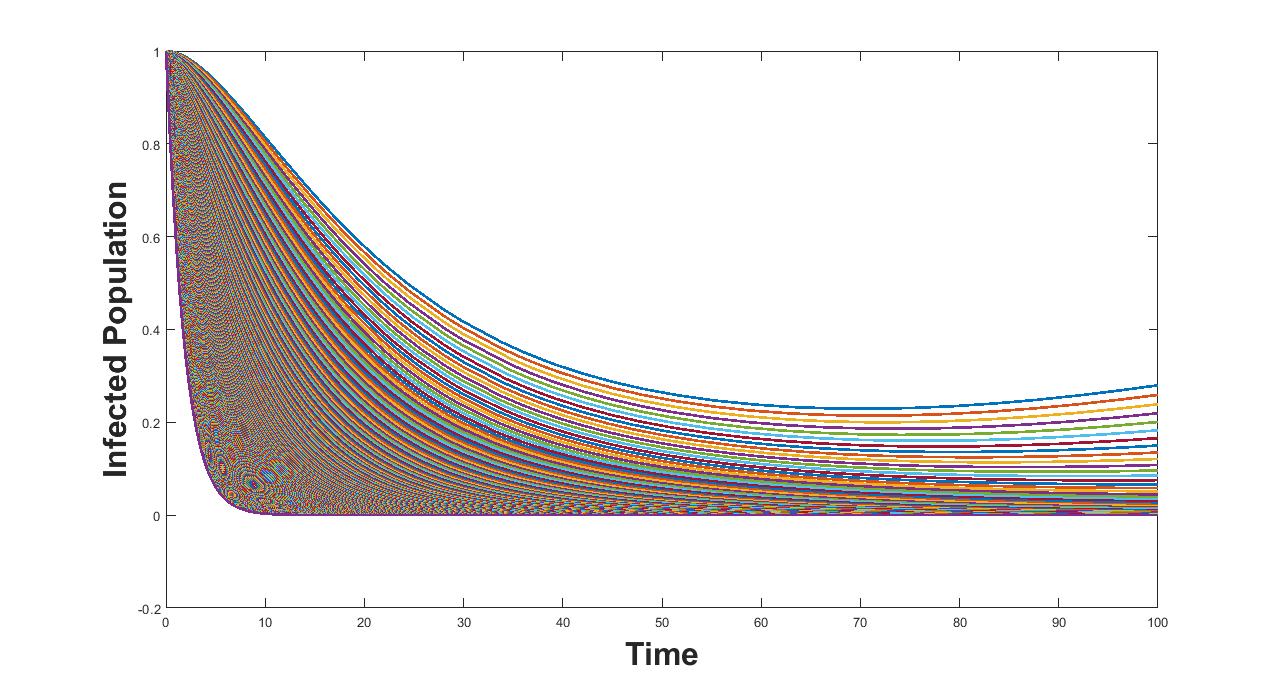}
				\hspace{-.4cm}
				\includegraphics[width=2in, height=1.8in, angle=0]{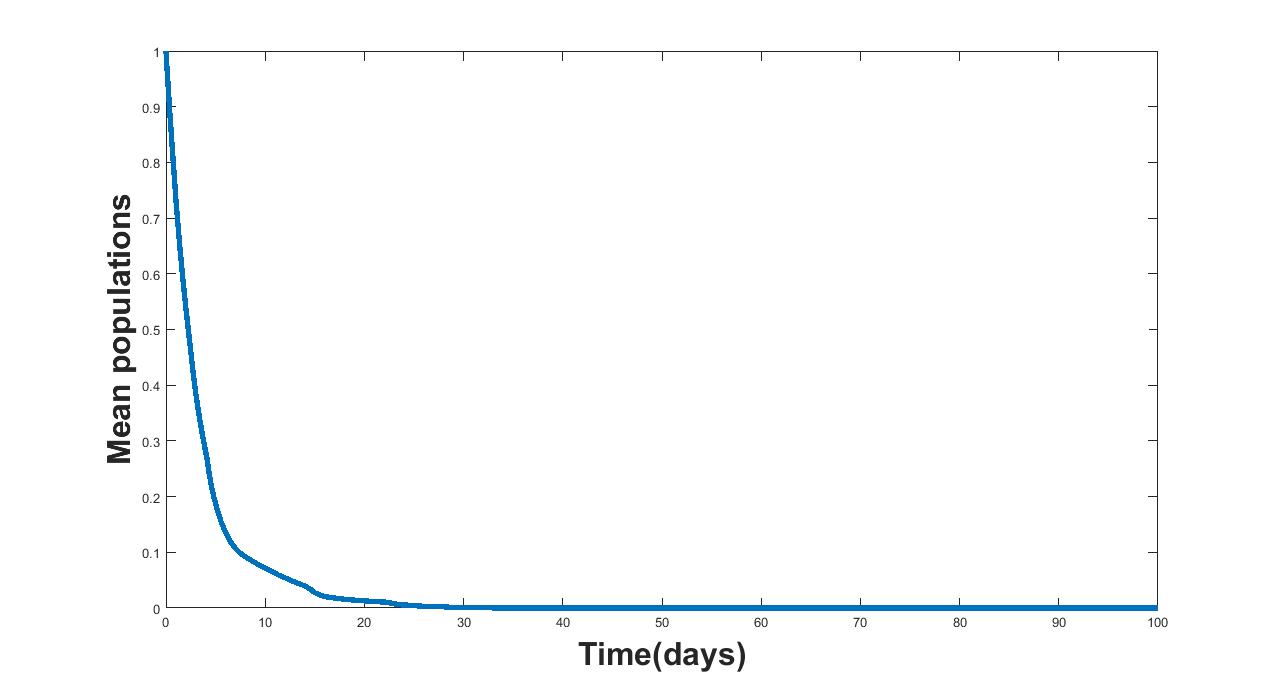}
				\hspace{-.395cm}
					\includegraphics[width=2in, height=1.8in, angle=0]{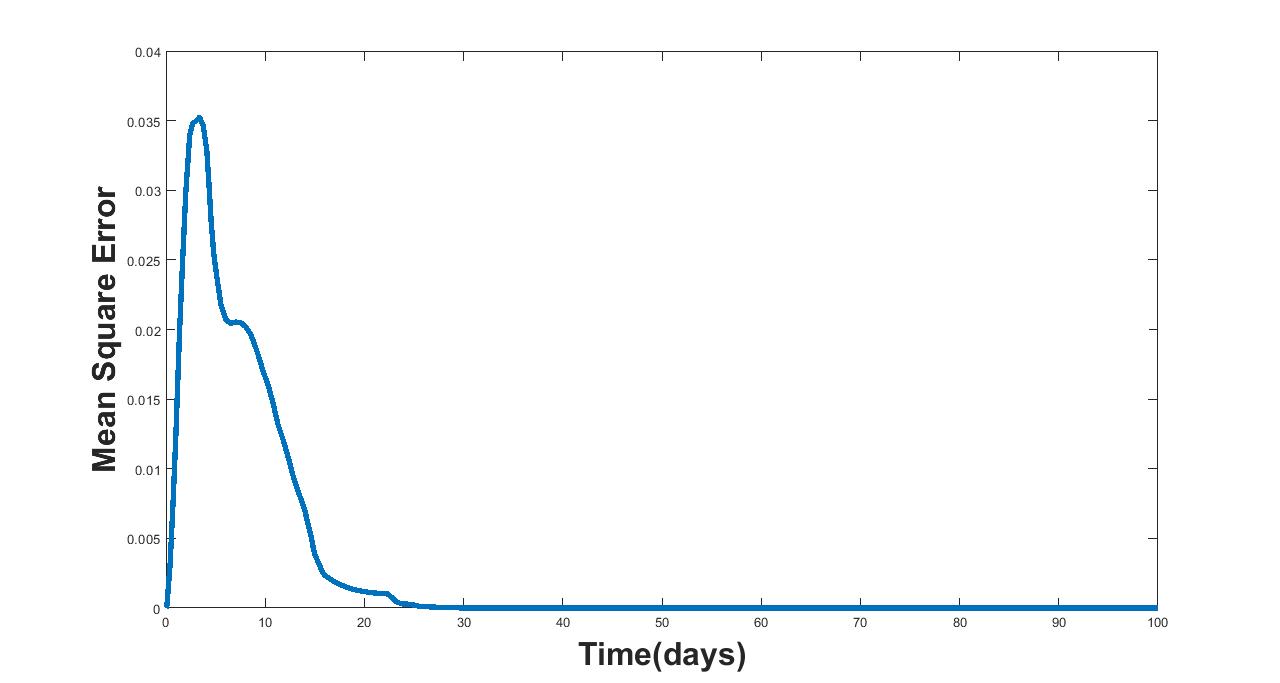}
				
			\caption*{(a) Interval I :  0 to 0.5}
				
			\end{center}
		\end{figure}

			\newpage
	\begin{table}[hbt!]
		\caption{Summary of Sensitivity Analysis}
		\centering
		\label{sen_anl}
		{
			\begin{tabular}{|l|l|l|}
				\hline
				\textbf{Parameter} & \textbf{Interval} & \textbf{Step Size}  \\
				\hline
				
				$\mu$ & 0 to 0.5 & \checkmark
					\\ \cline{2-3}
				 & .5 to 2.5  & $\times$\\
				 \hline
				 $\beta$ & 0 to 0.5 & \checkmark
					\\ \cline{2-3}
				 & 2  to 3  & $\times$
				 	\\
				
				 \hline
				 	$ \alpha $ & 0 to 1 & $\times$
					\\ \cline{2-3}
				 & 2 to 5  & $\times$ \\
				 \hline
				$\gamma$ & 0 to 1 & \checkmark
				\\ \cline{2-3} 
				& 1 to 2.5 & $\times$
				\\ \hline
				$\epsilon$ & 0 to 0.5 & \checkmark 
				\\ \cline{2-3}
				& 1.5 to 2.5 & $\times$ \\ 
				\hline 
					$\omega $ & 0 to 2 & \checkmark
				\\ \cline{2-3}
				& 2 to 4 &  \checkmark\\ 
				\hline 
					$\delta$ & 0 to 1 & $\times$
				\\ \cline{2-3}
				& 1 to 2.5 &  $\times$\\
					\hline 			
			\end{tabular}
		}
	\end{table}
	The summary of sensitivity analysis is given in Table \ref{sen_anl}. Due to brevity, the detailed sensitivity analysis of the rest of the parameters is given in Appendix 1.

\section{Numerical Simulations}
\label{sec5}
In this section with the help of the numerical simulations we verify our theoretical results. The chosen parameter values for the model are given in the the Table \ref{parameters table}. We have used the MATLAB software for the simulation.

\begin{figure}[hbt!]
\centering
\includegraphics[width=5in, height=3in, angle=0]{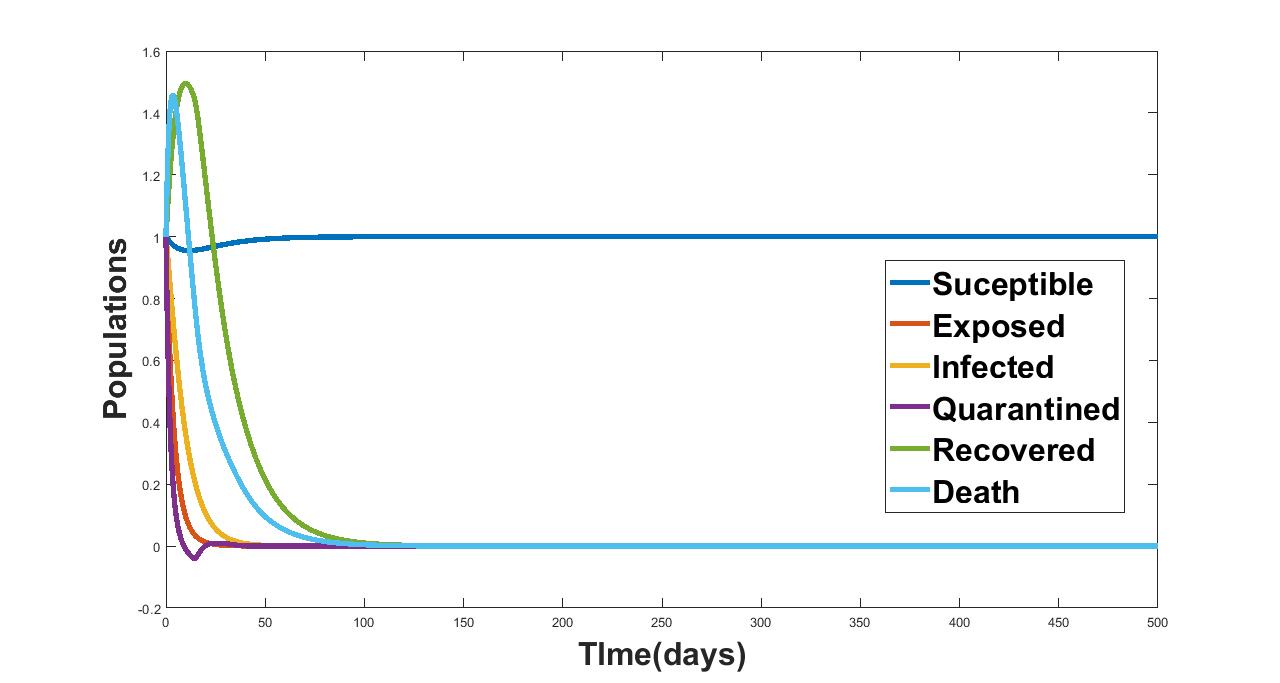}
\caption{Stability of $E_0$ for $\tau^*=5$}
\label{f1}
\end{figure}

\begin{figure}[hbt!]
\centering
\includegraphics[width=5in, height=3in, angle=0]{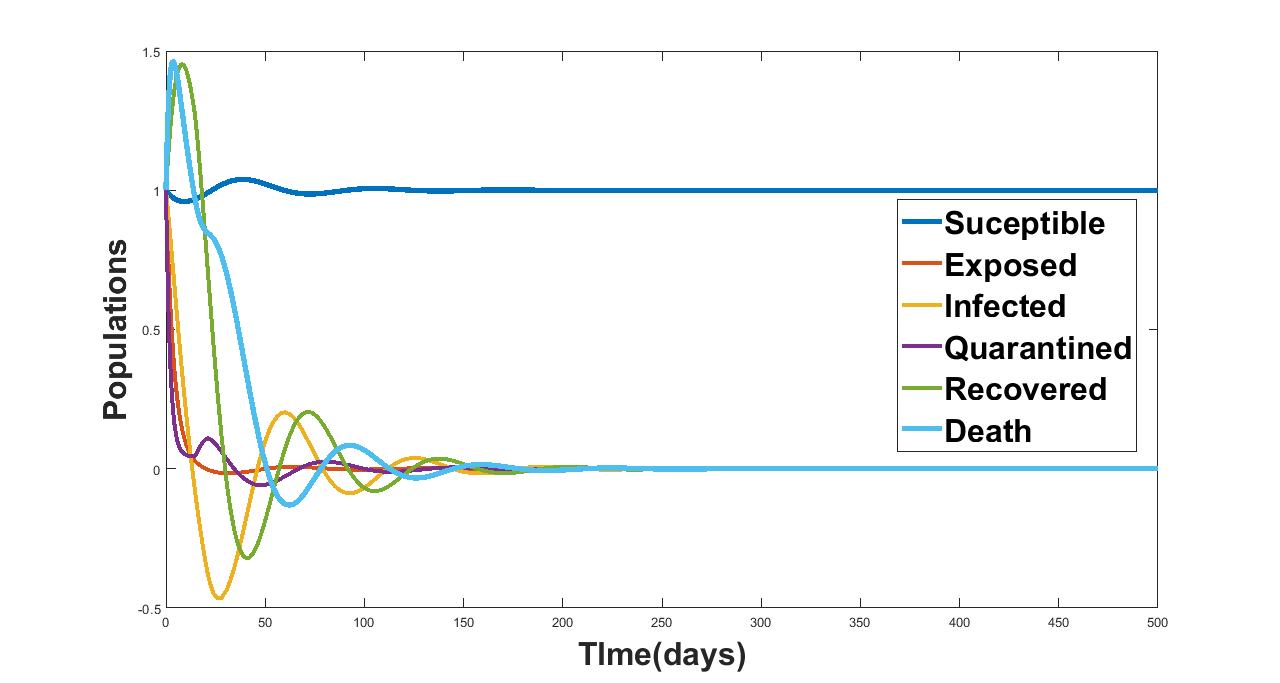}
\caption{Stability of $E_0$ for $\tau^*=20$}
\label{f2}
\end{figure}

For the parameter values from Table \ref{parameters table}, we find that $f_1>0$ and $(\mu+\delta \alpha)^2 > \rho^2 (1-\alpha)^2$. The values of $d_2^2$ was found to be 0026 and $e_2^2$ was found to be 0.00045. Therefore the condition $d_2^2 > e_2^2$ of Theorem \ref{theorem} is satisfied. Now from Theorem \ref{theorem}, $E_0$ should be asymptotically stable for all the possible values of delay. Here we take two values of delay and check for the stability of $E_0$. From Figure \ref{f1} , Figure \ref{f2}, we see that  the infection free equilibrium $E_0$ is locally asymptotically stable.

\newpage

\begin{figure}[hbt!]
\centering
\includegraphics[width=5in, height=3in, angle=0]{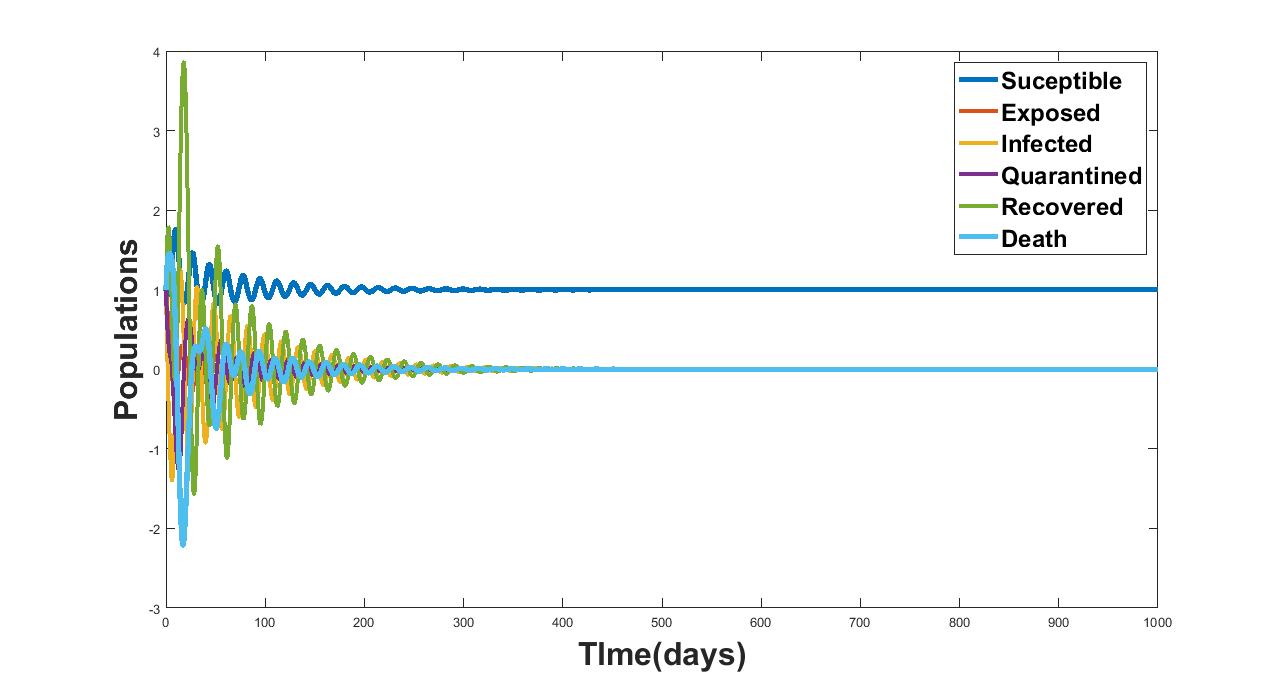}
\caption{Stability of $E_0$ for $\tau = 5 < \tau^*=5.8485$}
\label{f3}
\end{figure}

\begin{figure}[hbt!]
\centering
\includegraphics[width=5in, height=3in, angle=0]{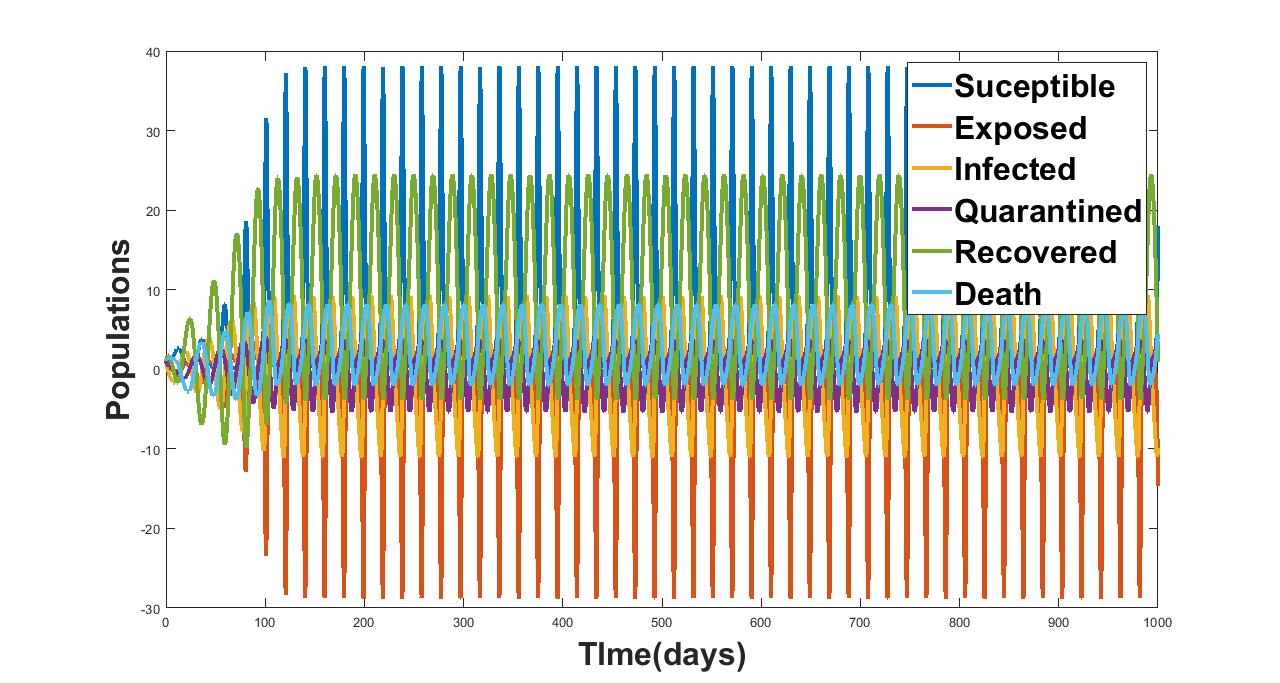}
\caption{Stability of $E_0$ for $\tau =7 > \tau^*=5.8485$}
\label{f4}
\end{figure}

\begin{figure}[hbt!]
\centering
\includegraphics[width=5in, height=3in, angle=0]{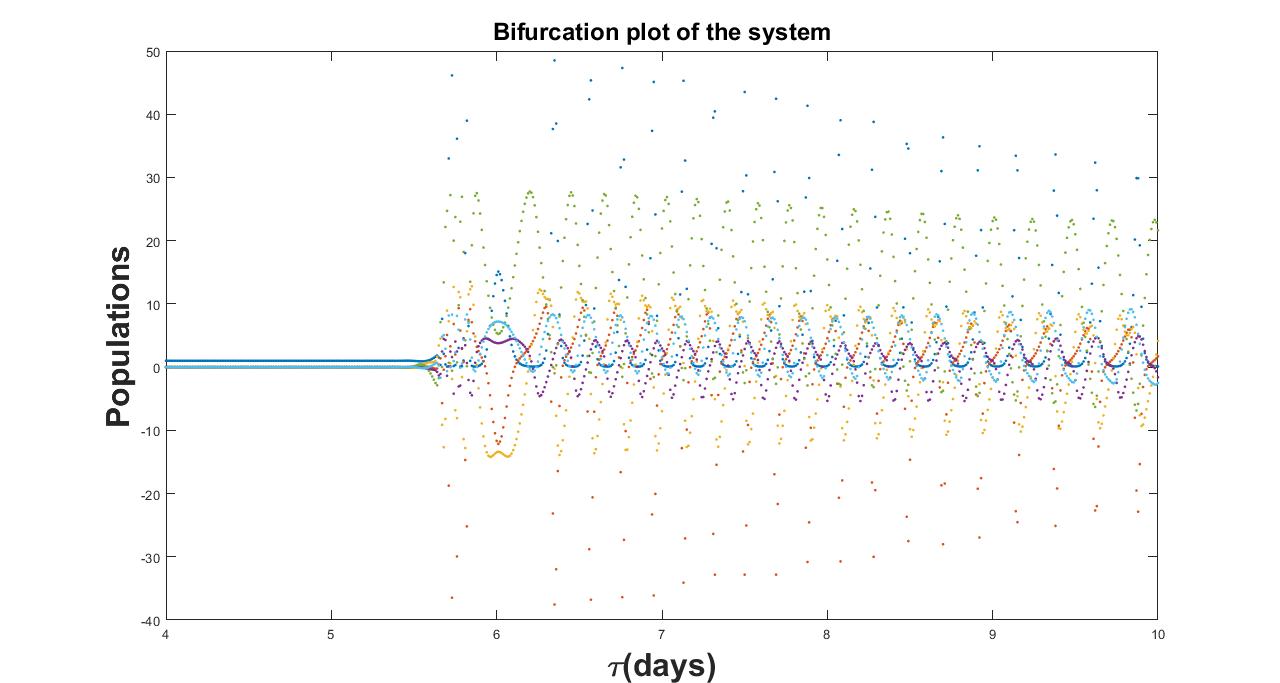}
\caption{Stability of $E_0$ for $\tau^*=5.8485$}
\label{f5}
\end{figure}

With the parameter values from Table \ref{parameters table}, we have  $d_2^2= 0.0062$ and $e_2^2= 0.0069$. Therefore we have $d_2^2 < e_2^2$. Also $\gamma^2 =  0.0142 > \omega^2=0.0139$ and $d_2 =0.0128 < \omega^2=0.0139$. We also see that the condition $x +y= 9.4271 > z =9.3946$ is satisfied. The critical value of delay $\tau^*$ is calculated to be  5.8485. Therefore from Theorem \ref{theorem}, we have $E_0$ to be asymptotically stable for $\tau<\tau^*$ and unstable for $\tau>\tau^*$. From Figure \ref{f3} , Figure \ref{f4}, we see that $E_0$ is locally asymptotically stable for $\tau=5 < \tau^*=5.8485$ and unstable for $\tau=7 > \tau^* = 5.8485$. We plot the bifurcation graph in Figure \ref{f5}. From Figure \ref{f5}, we see that the solutions of the system changes its nature from being asymptotically stable to unstable as the value of $\tau$ crosses 5.5. Therefore Hopf bifurcation occurs at $\tau^*=5.8485$.

\newpage
     
\begin{figure}[hbt!]
\centering
\includegraphics[width=5in, height=3in, angle=0]{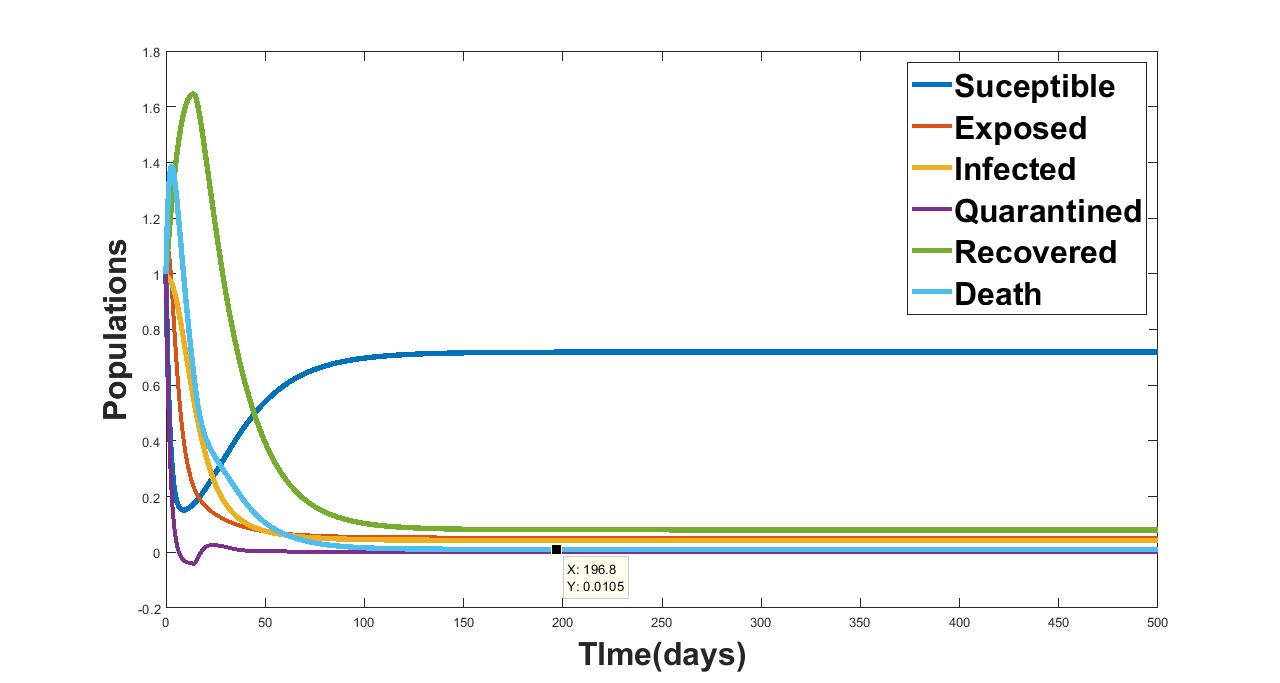}
\caption{Stability of $E^*$ for $\tau^*=10$}
\label{f6}
\end{figure}

\begin{figure}[hbt!]
\centering
\includegraphics[width=5in, height=3in, angle=0]{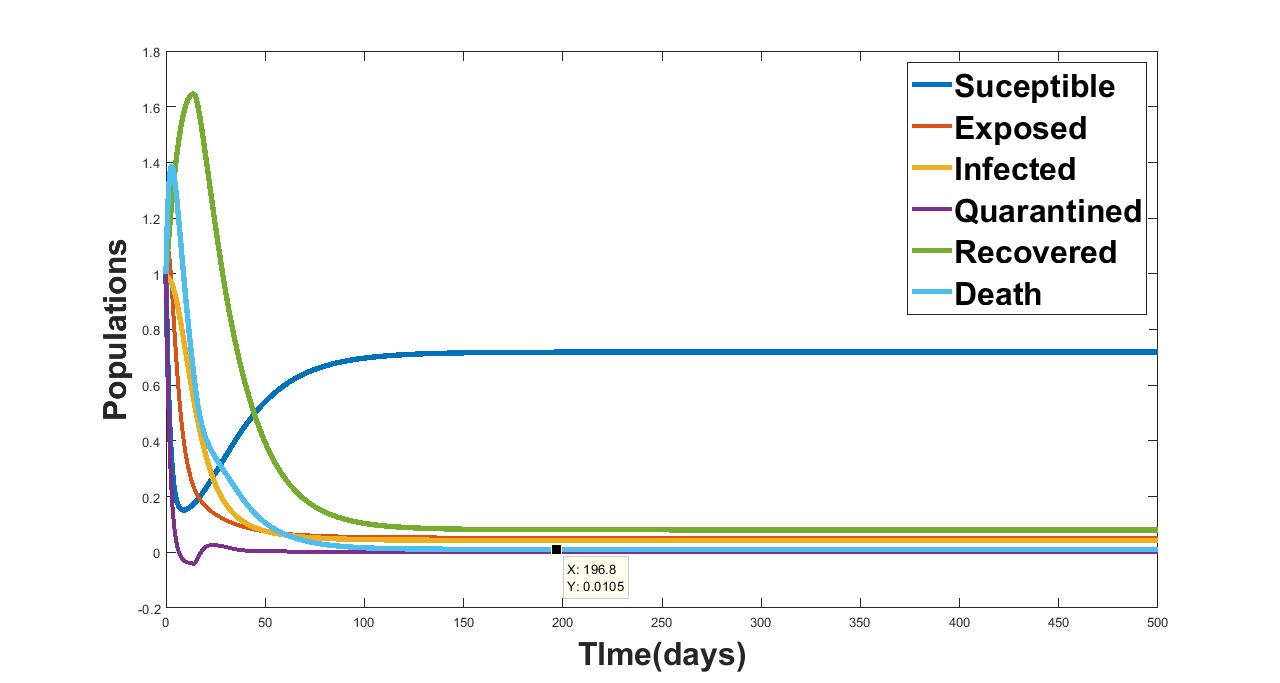}
\caption{Stability of $E^*$ for $\tau^*=20$}
\label{f7}
\end{figure}

Now we numerically check the stability of infected equilibrium $E^*$. For the parameter values from Table \ref{parameters table}, the infected equilibrium $E^*$ was found to be $(0.633,0.093,0.08,0.003,0.16,0.02)$ and we calculate the values of $\beta \epsilon = .0882$ and $(\epsilon+\mu)*(\gamma+\mu+\rho*exp^{-\gamma*\tau})=.0559$. Therefore  the existence of $E^*$ is guaranteed. The values of $c_1, c_2, c_3$ were calculated to be  $0.22,0.000035, 0.0000025$ which are all positive. Also we calculate the values of $(\mu + \delta*\alpha)^2$ to be 0.247 and value of $\rho^2*(1-\alpha)^2$ to be 0.0017 and we see that $(\mu + \delta*\alpha)^2 > \rho^2*(1-\alpha)$. Therefore from the stability analysis of $E^*$ we must have $E^*$ to be locally asymptotically stable for all delays. From Figure \ref{f6} and Figure \ref{f7}, we see that $E^*$ is stable for all delays.\\

\begin{figure}[hbt!]
\centering
\includegraphics[width=5in, height=3in, angle=0]{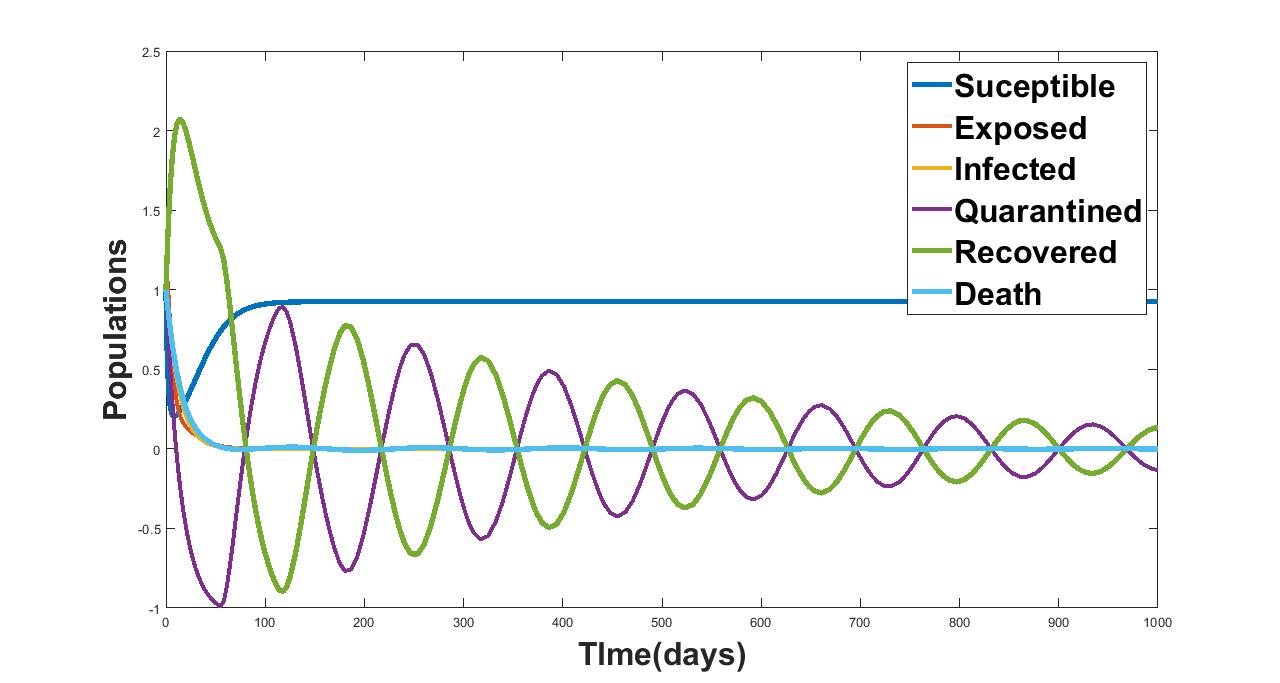}
\caption{Stability of $E^*$ for $\tau^*=20$}
\label{f8}
\end{figure}

When we take the value of $\alpha$ to be 0.001 we see that $(\mu + \delta*\alpha)^2 = 0.004 < \rho^2*(1-\alpha)=0.005$. This violates the stability condition for $E^*$ and we expect $E^*$ to be unstable. From Figure \ref{f8}, we see that $E^*$ starts loosing its stability and starts becoming unstable when we have $(\mu + \delta*\alpha)^2 < \rho^2*(1-\alpha)$.

\newpage

\section{Effect of Temperature}
\label{sec6}
\subsection{Impact of Temperature : }

\begin{figure}[hbt!]
\centering
\includegraphics[width=5in, height=3in, angle=0]{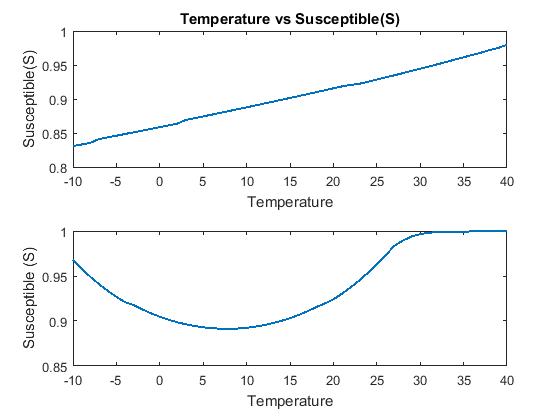}
\caption{Susceptible vs temperature}\label{ts}
\end{figure}

\begin{figure}[hbt!]
\centering
\includegraphics[width=5in, height=3in, angle=0]{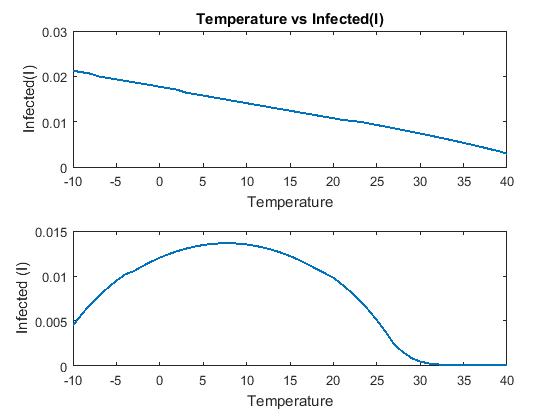}
\caption{Infected vs temperature}\label{ti}
\end{figure}

\begin{figure}[hbt!]
\centering
\includegraphics[width=5in, height=3in, angle=0]{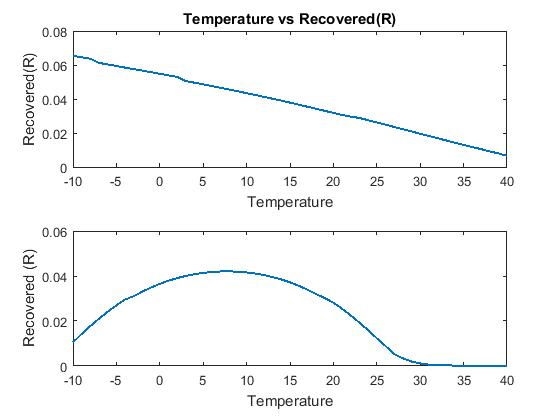}
\caption{Recovered vs temperature}\label{tr}
\end{figure}

\begin{figure}[hbt!]
\centering
\includegraphics[width=5in, height=3in, angle=0]{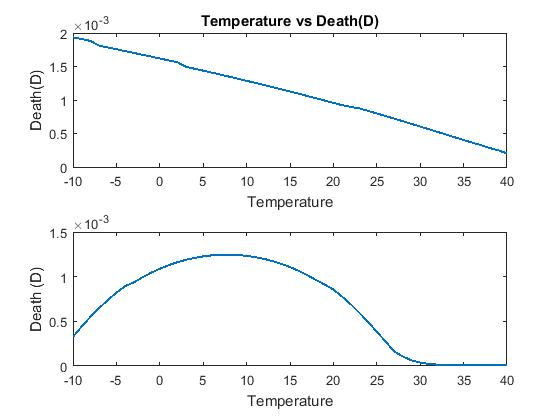}
\caption{Death vs temperature}\label{td}
\end{figure}

Temperature is one of the key indicators that influence the spread of SARS-CoV-2 virus. For temperatures ranging from -10 to 40, the average population belonging to each compartment is calculated using dde23 method in MATLAB. Since we represented the spread parameter $\beta$ as a linear and a quadratic function of temperature, the first subplot gives the average values with respect to $\beta$ as a linear function of temperature. And the second subplot gives the average values with respect to $\beta$ as a quadratic function of temperature. Both the subplots in Figures \ref{ts}, \ref{ti}, \ref{tr}, \ref{td} concluded that with the decrease in temperature, the average number of infections and the average number of deaths in a population of individuals increases. 

\newpage
\newpage

\subsection{Impact of isolation : }

\begin{figure}[hbt!]
\centering
\includegraphics[width=5in, height=3in, angle=0]{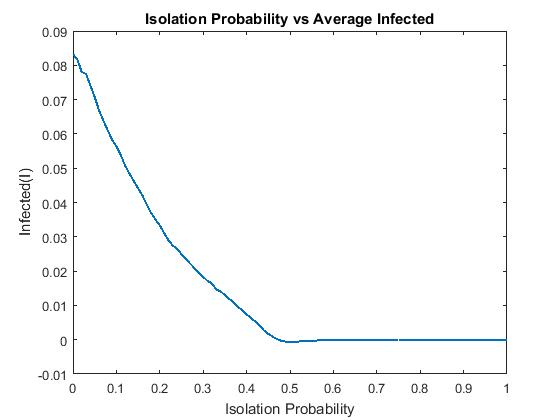}
\caption{Isolation probability vs Infected}\label{pi}
\end{figure}

\begin{figure}[hbt!]
\centering
\includegraphics[width=5in, height=3in, angle=0]{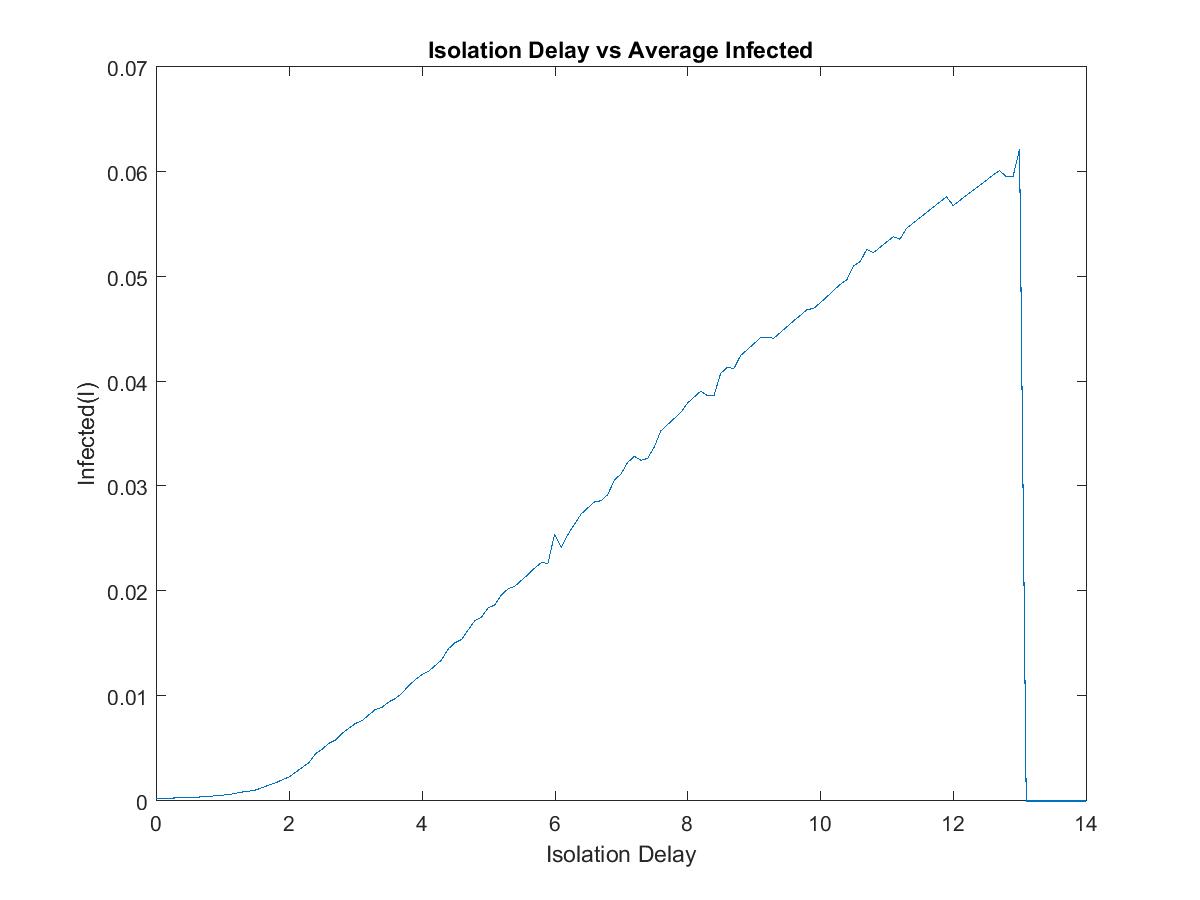}
\caption{Isolation delay vs Infected}\label{taui}
\end{figure}

Isolation of infected individuals plays a very crucial role in controlling the pandemic situations. With the infections highly probable to increase in the coming winter season, efficient isolation strategies play an important role in controlling the spike of new infections. In addition to the probability of isolation, the delay in isolating infected individuals also play very crucial role in controlling infections. Figure \ref{pi} and Figure \ref{taui} give an idea on the impact of isolation probability and isolation delay with respect to the new infections.

\section{Conclusions}
\label{sec7}
The system of delay differential equations proposed in this paper, incorporates key characteristics of COVID-19. This includes the impact of temperature on the spread parameter $\beta$ and the delay in isolating infected individuals. The reproduction number is calculated and the impact of temperature and two isolation indicators(isolation probability and isolation delay) with the temperature are plotted. The proposed model admits two equilibria, i.e., disease free equilibrium and infection equilibrium. The conditions for the stability of equilibrium points is obtained. The system exhibits hopf bifurcation. The critical delay at which hopf bifurcation occurs is calculated. Sensitivity analysis is performed for the proposed model. It is identified that a few parameters are sensitive only at some interval. The impact of temperature and isolation on COVID-19 is observed that with the decrease in temperature, the average number of infections and the average number of deaths increases. This observation is in similar lines with the conclusions in \cite{temp1} that 60.0\% of the confirmed cases of coronavirus disease 2019 (COVID-19) occurred in places where the air temperature ranged from 5 C to 15 C, with a peak in cases at 11.54 C. With the winter season round the corner, these results warn us of an increase in daily new infected cases and possibly a second wave of COVID-19. Hence, the measures of face masks, sanitizers and social distancing has to be strictly followed in the coming future to avoid any sudden spike of cases.
 
\newpage

\bibliographystyle{elsarticle-num}
\bibliography{mybibfile}

\newpage
\section*{Appendix 1: Sensitivity Analysis for other parameters}
\label{sec9}

\subsection{{$\alpha$}}
		
\begin{figure}[hbt!]
\begin{center}
\includegraphics[width=2in, height=1.8in, angle=0]{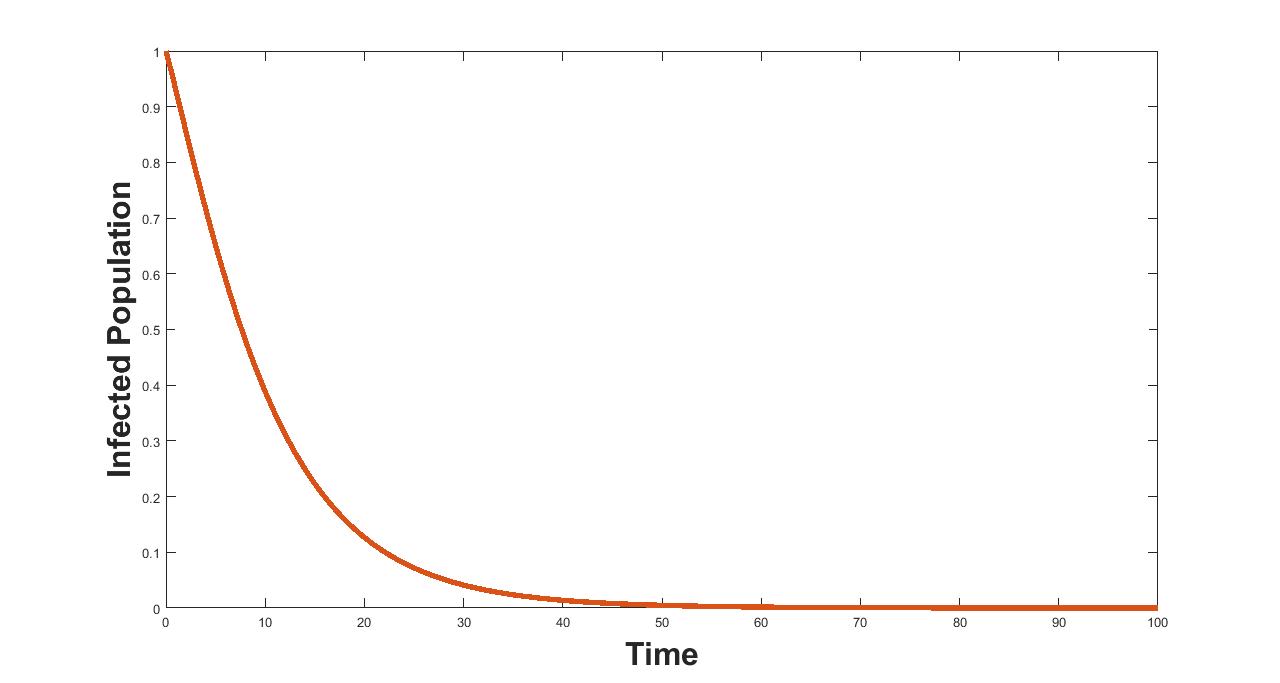}
\hspace{-.4cm}
\includegraphics[width=2in, height=1.8in, angle=0]{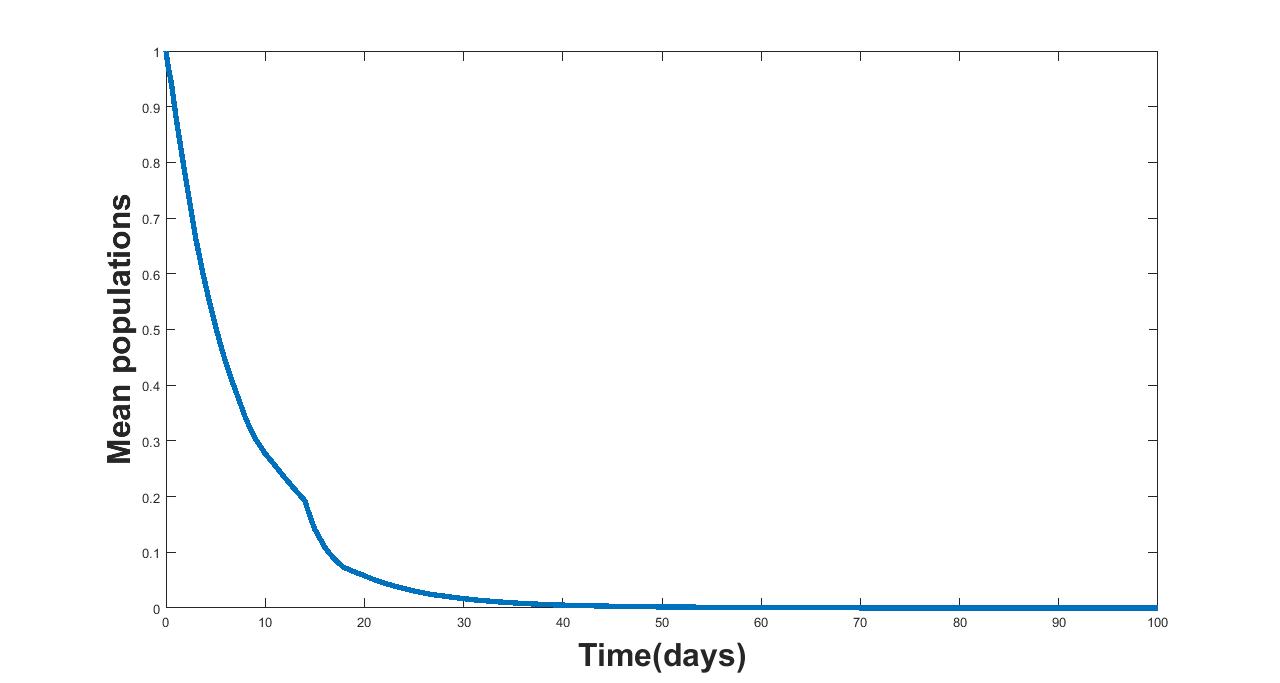}
\hspace{-.395cm}
\includegraphics[width=2in, height=1.8in, angle=0]{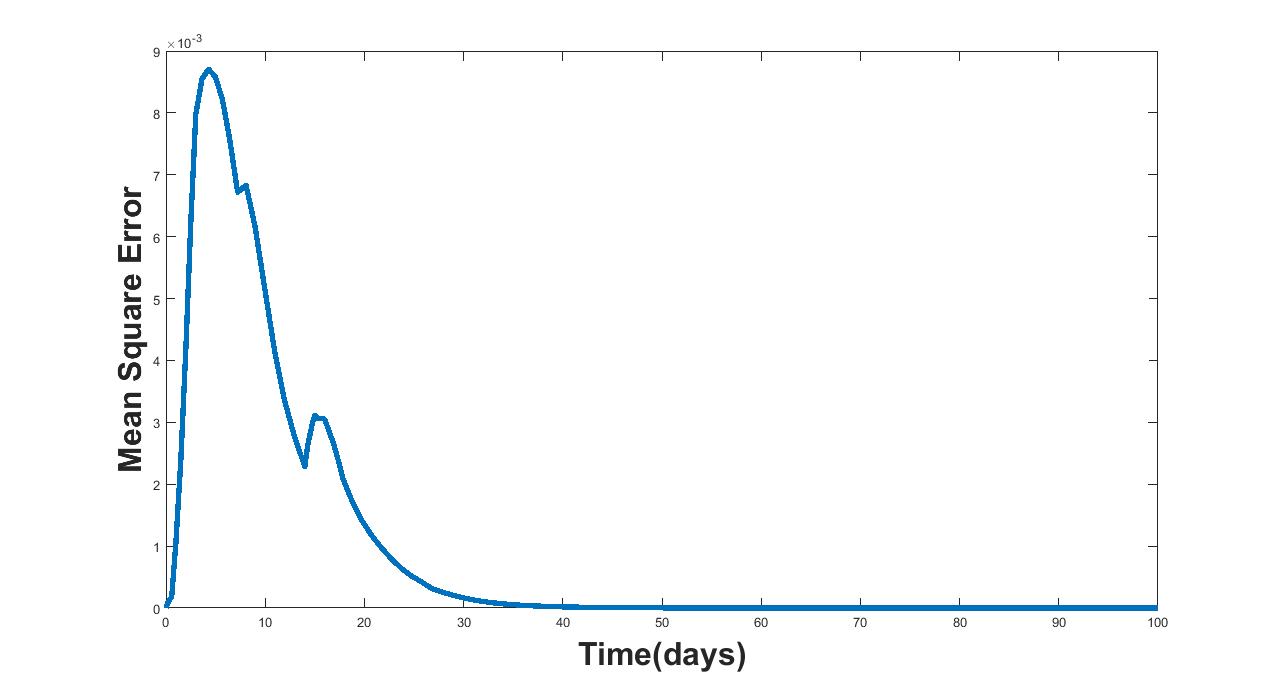}
\caption*{(a) Interval I : 0 to 3}
\end{center}
\end{figure}

\vspace{-3mm}

\begin{figure}[hbt!]
\begin{center}
\includegraphics[width=2in, height=1.8in, angle=0]{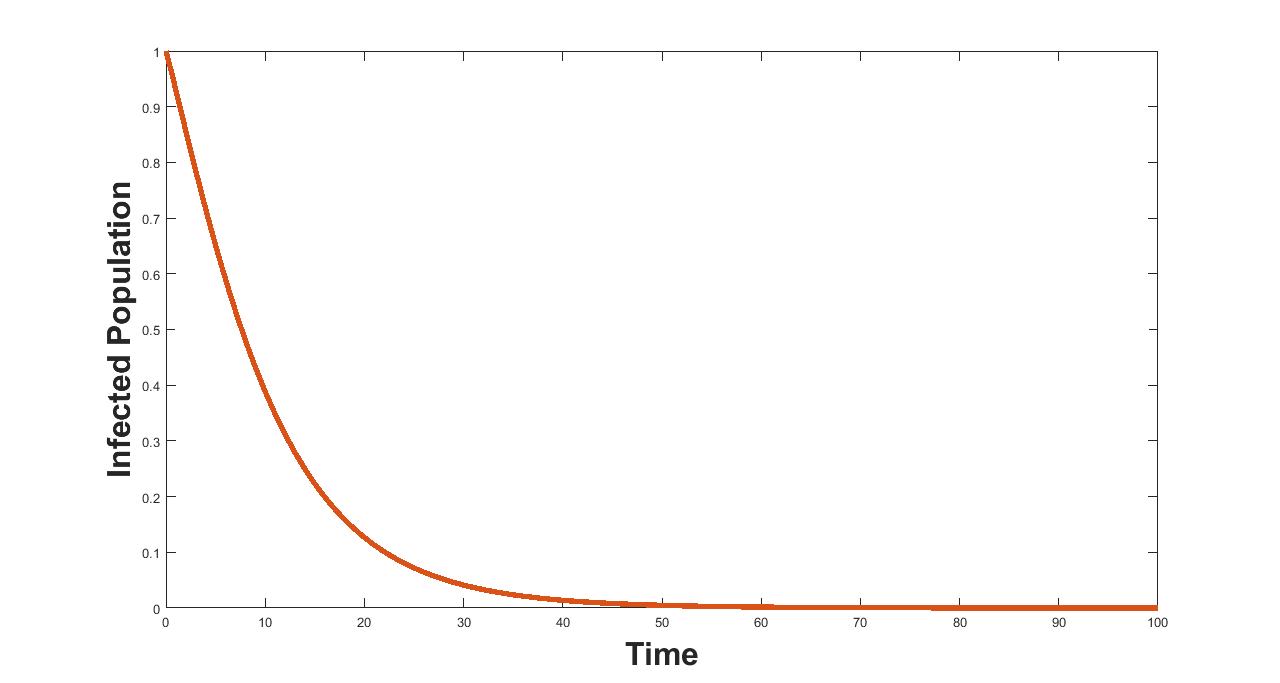}
\hspace{-.4cm}
\includegraphics[width=2in, height=1.8in, angle=0]{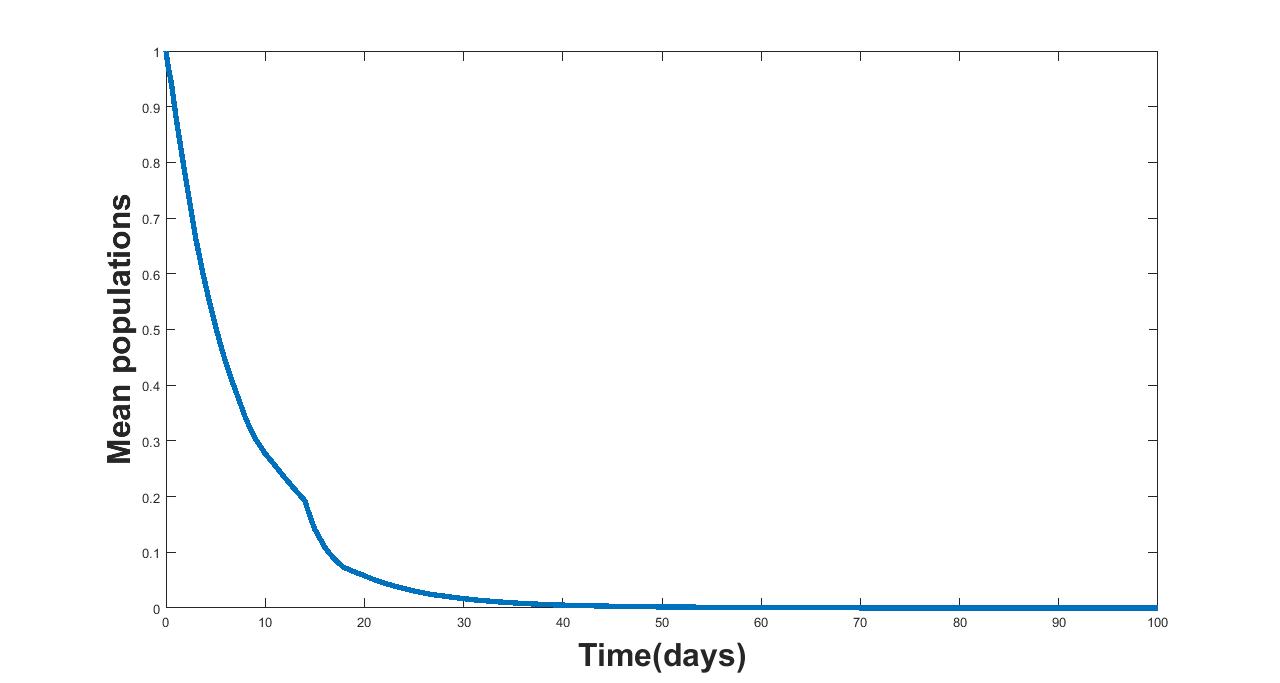}
\hspace{-.395cm}
\includegraphics[width=2in, height=1.8in, angle=0]{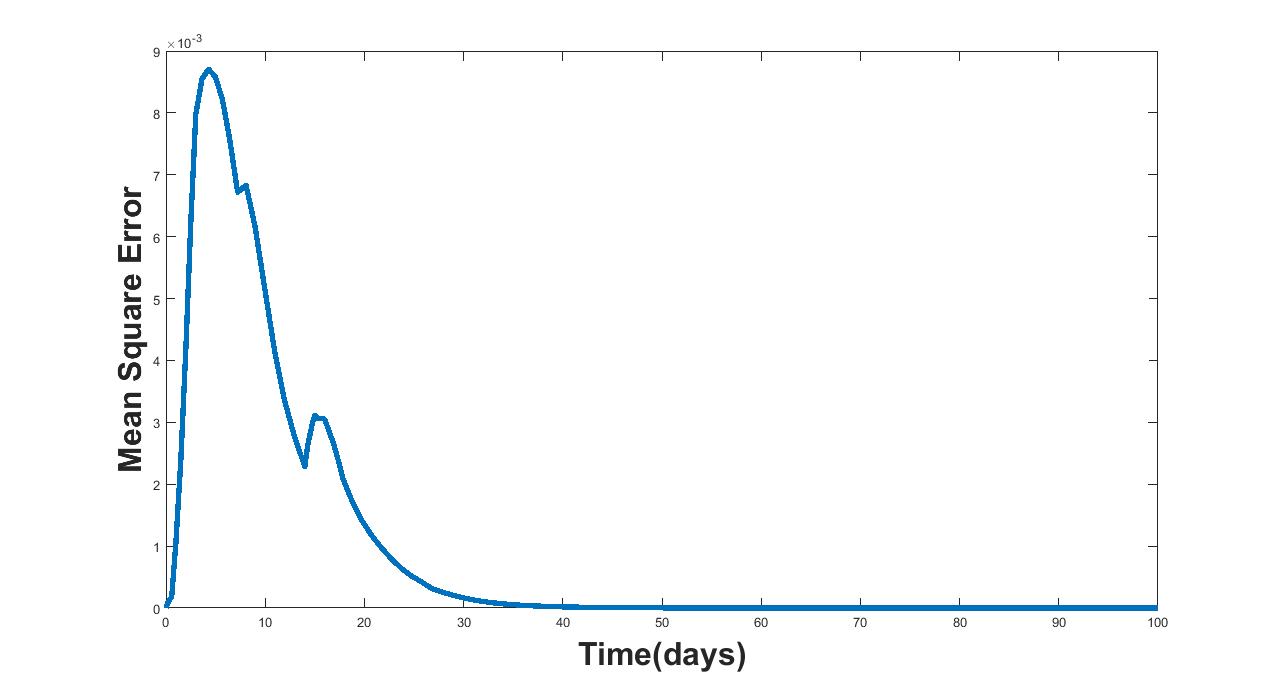}
\caption*{(b) Interval II : 3 to 5}

\vspace{5mm}
\caption{Sensitivity Analysis of $\alpha$}
\label{sen_alpha}
\end{center}
\end{figure}

\newpage

\subsection{Parameter $\boldsymbol{\epsilon}$}
		
\begin{figure}[hbt!]
\begin{center}
\includegraphics[width=2in, height=1.8in, angle=0]{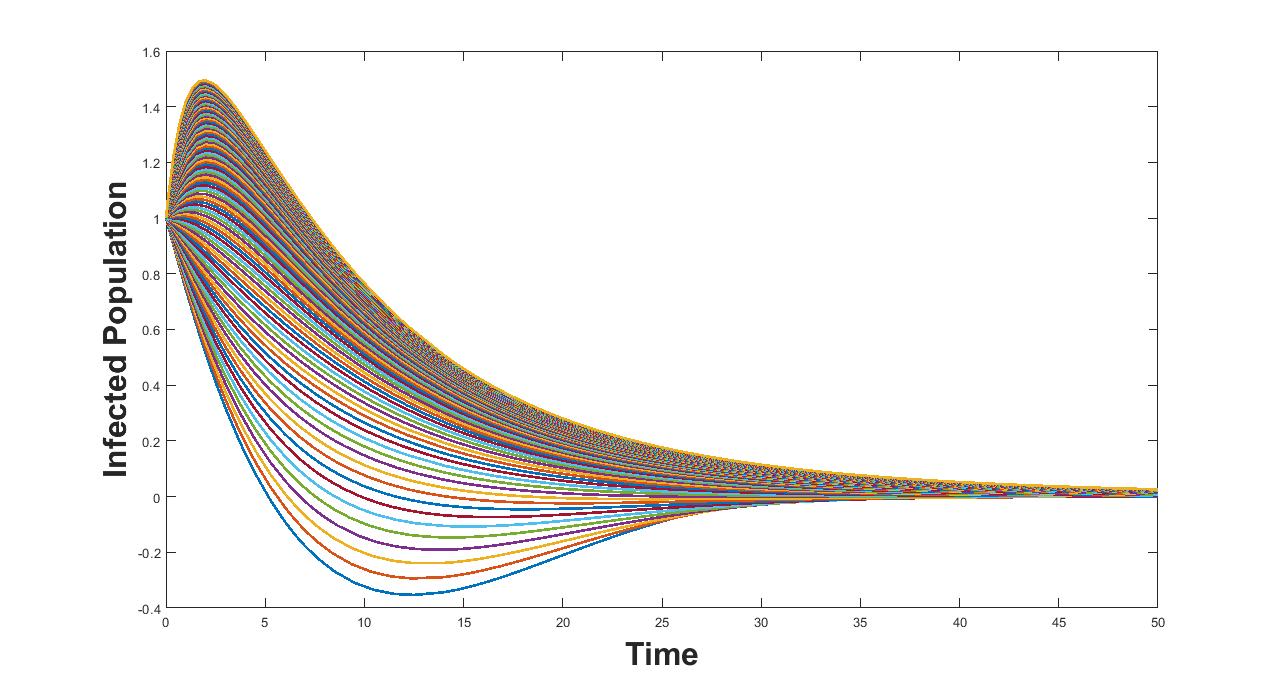}
\hspace{-.4cm}
\includegraphics[width=2in, height=1.8in, angle=0]{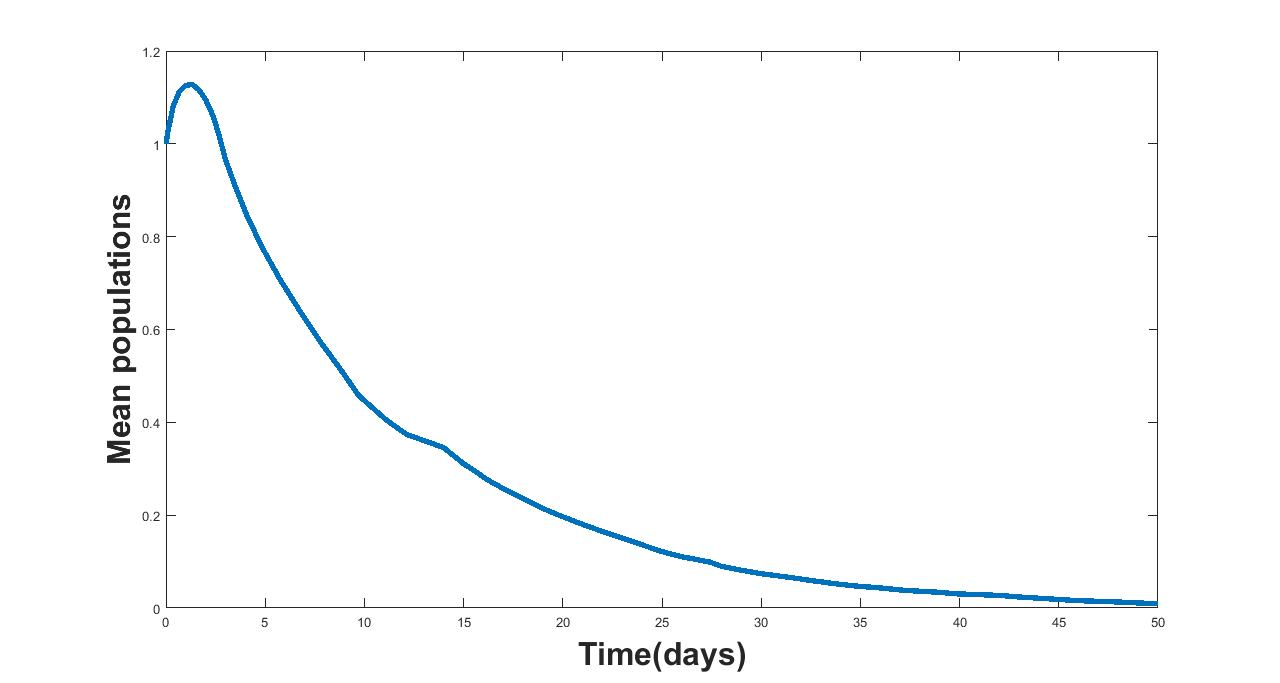}
\hspace{-.395cm}
\includegraphics[width=2in, height=1.8in, angle=0]{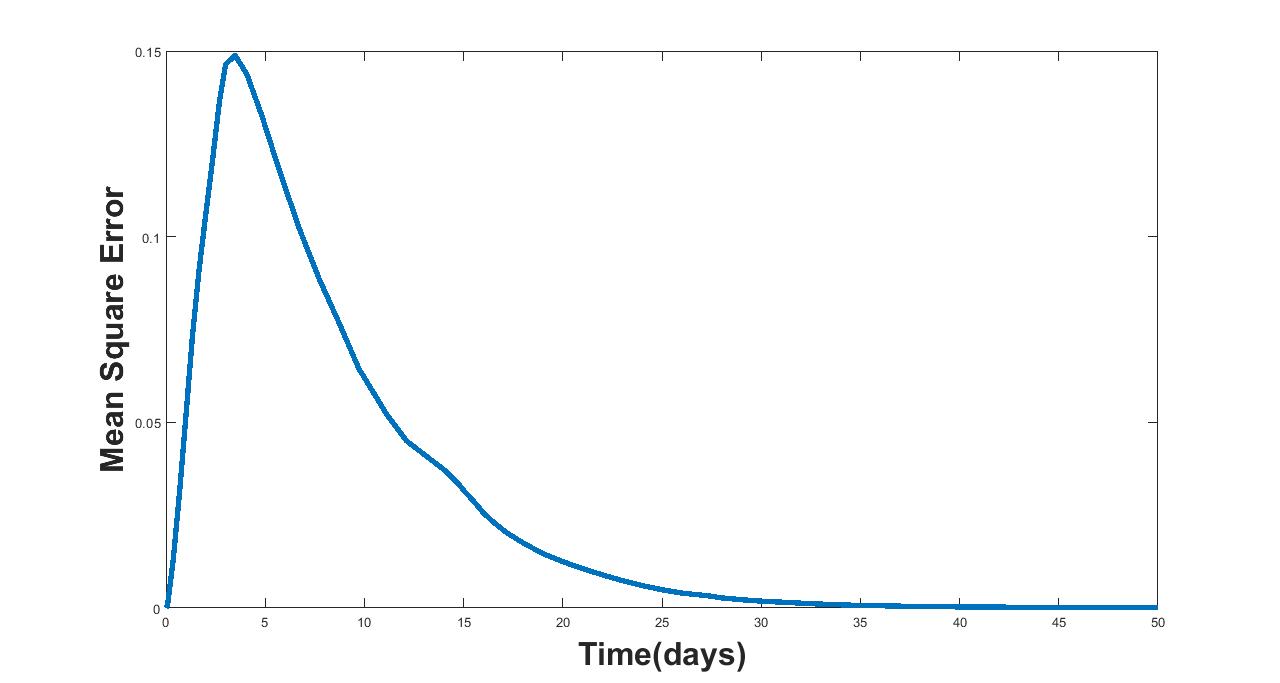}
\caption*{(a) Interval I : 0 to 0.5}
\end{center}
\end{figure}

\vspace{-3mm}

\begin{figure}[hbt!]
\begin{center}
\includegraphics[width=2in, height=1.8in, angle=0]{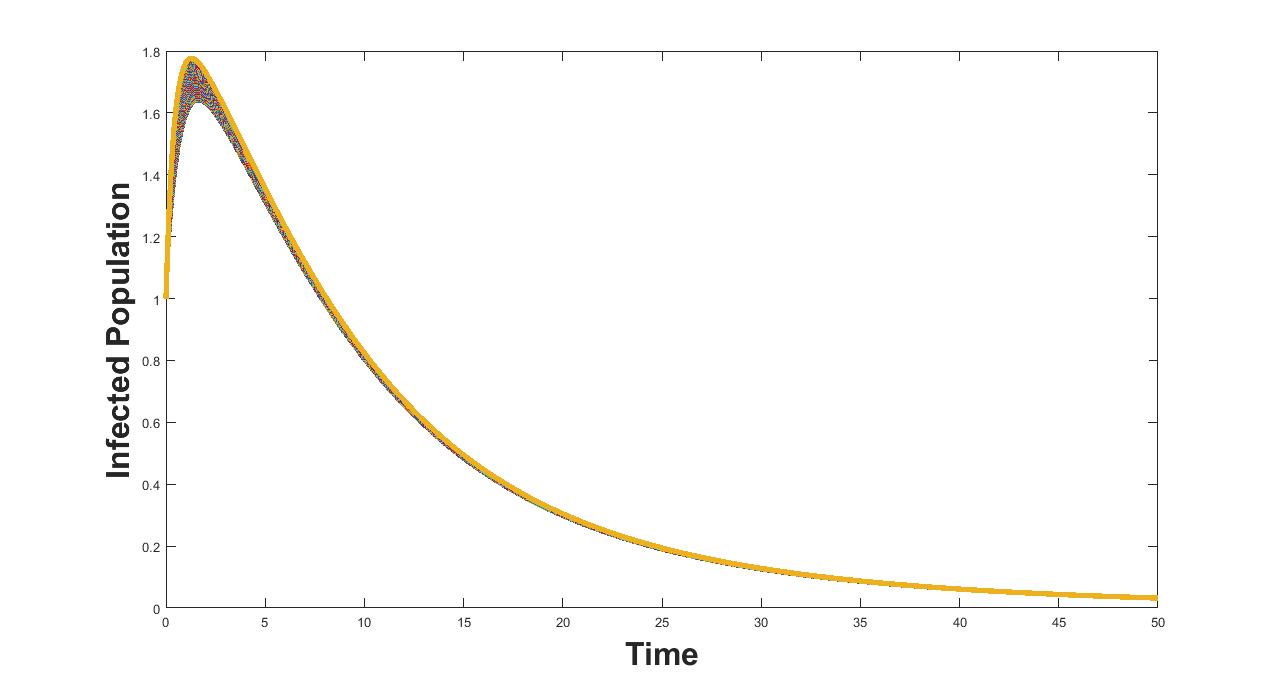}
\hspace{-.4cm}
\includegraphics[width=2in, height=1.8in, angle=0]{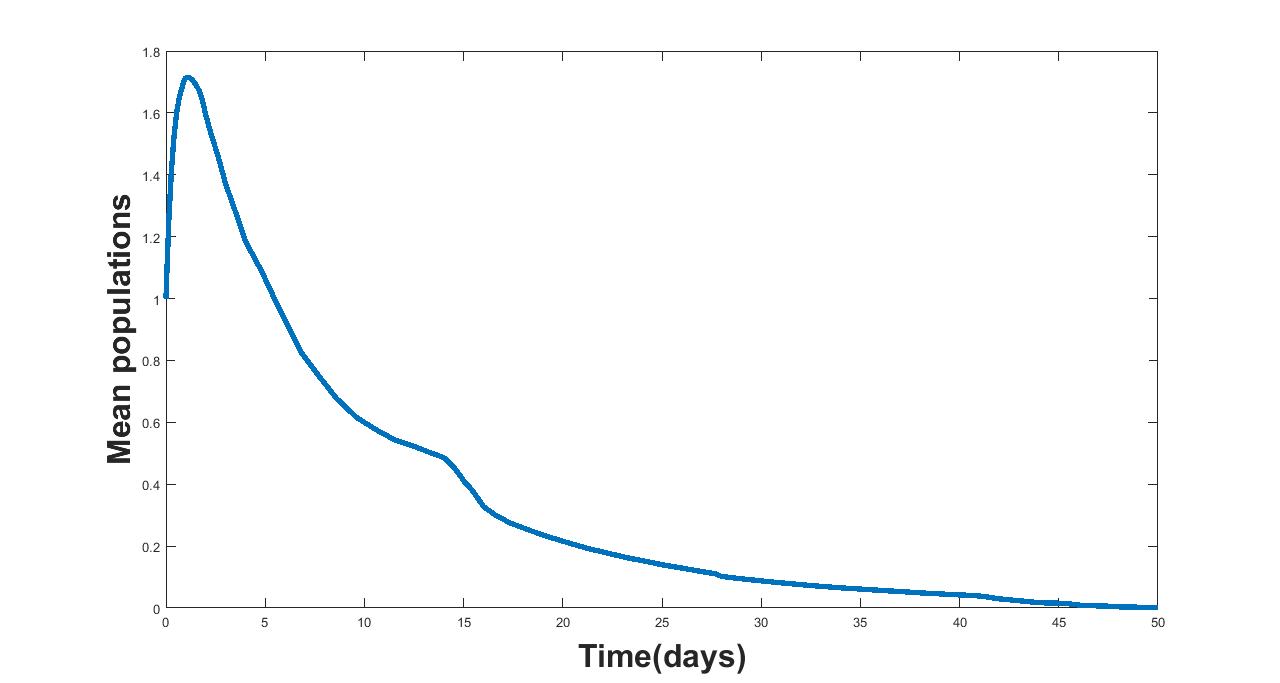}
\hspace{-.395cm}
\includegraphics[width=2in, height=1.8in, angle=0]{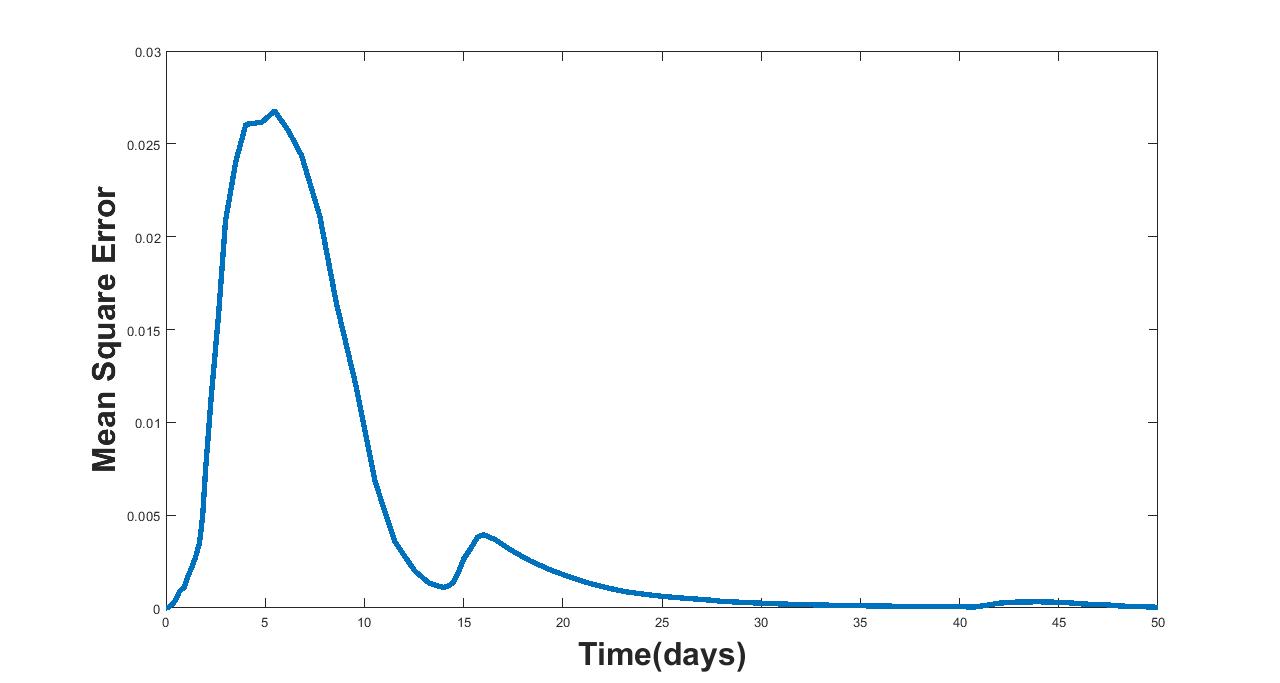}
\caption*{(b) Interval I : 1.5 to 2.5}

\vspace{5mm}
\caption{Sensitivity Analysis of $\epsilon$}
\label{sen_eps}
\end{center}
\end{figure}

\newpage

\subsection{Parameter $\boldsymbol{\delta}$}
		
\begin{figure}[hbt!]
\begin{center}
\includegraphics[width=2in, height=1.8in, angle=0]{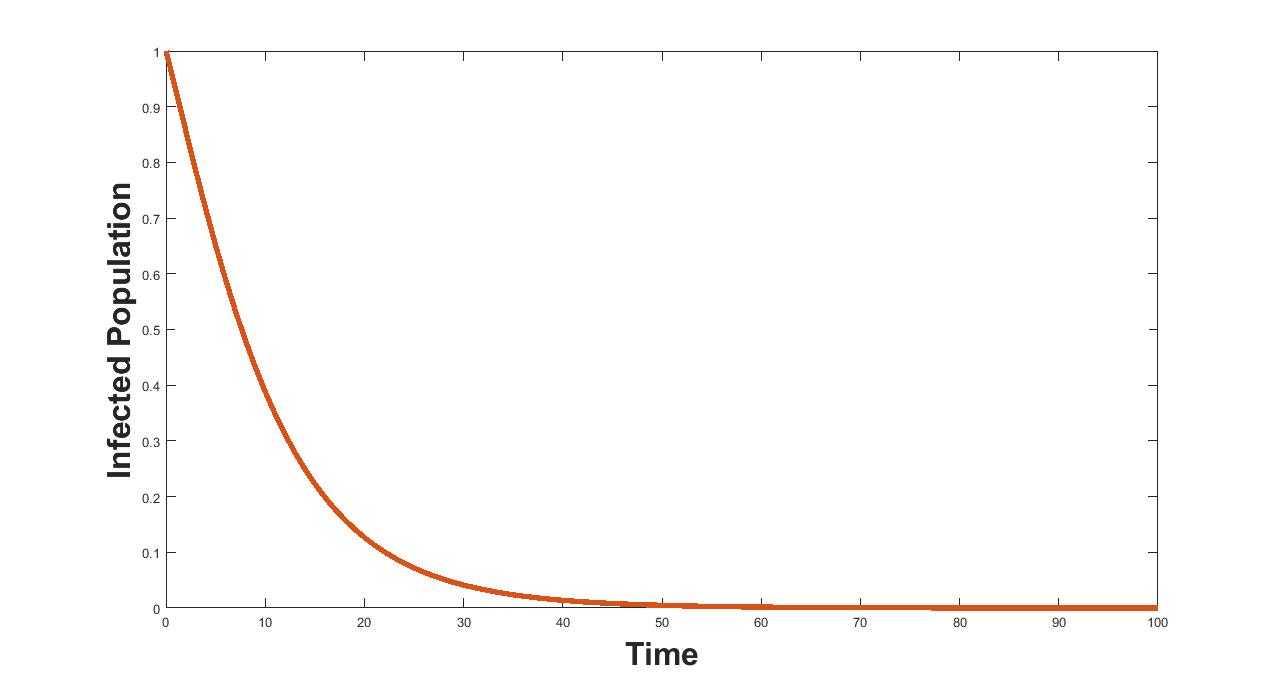}
\hspace{-.4cm}
\includegraphics[width=2in, height=1.8in, angle=0]{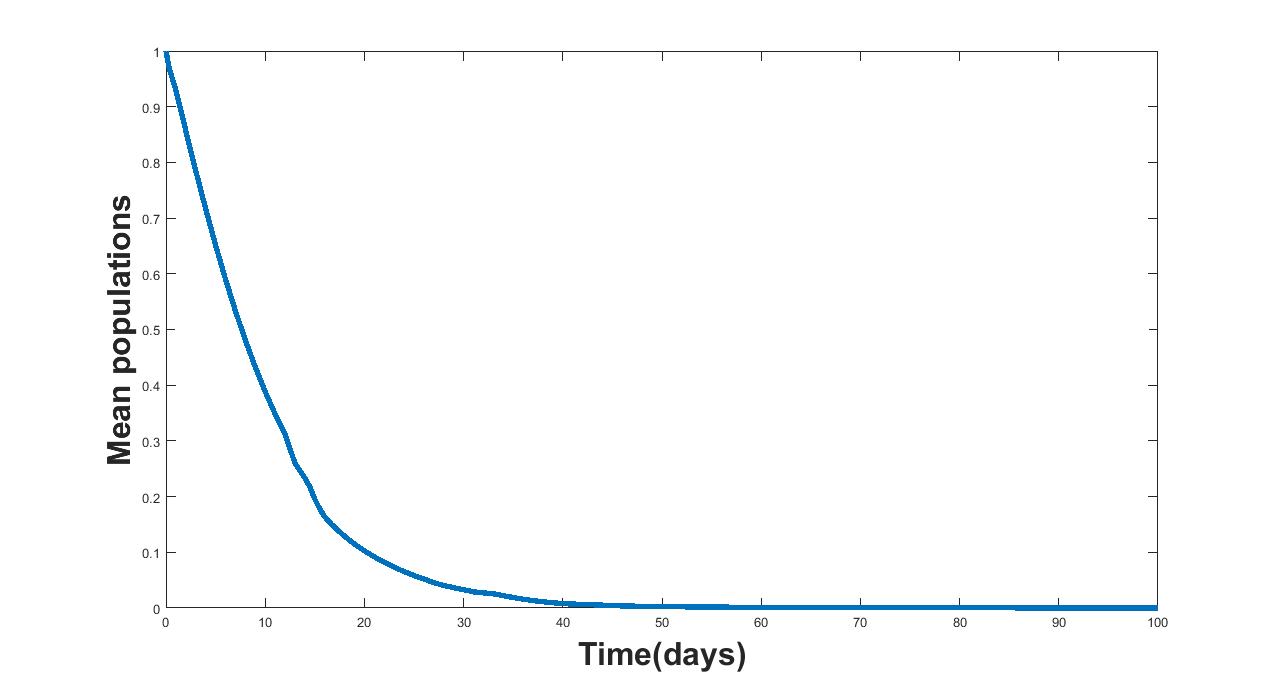}
\hspace{-.395cm}
\includegraphics[width=2in, height=1.8in, angle=0]{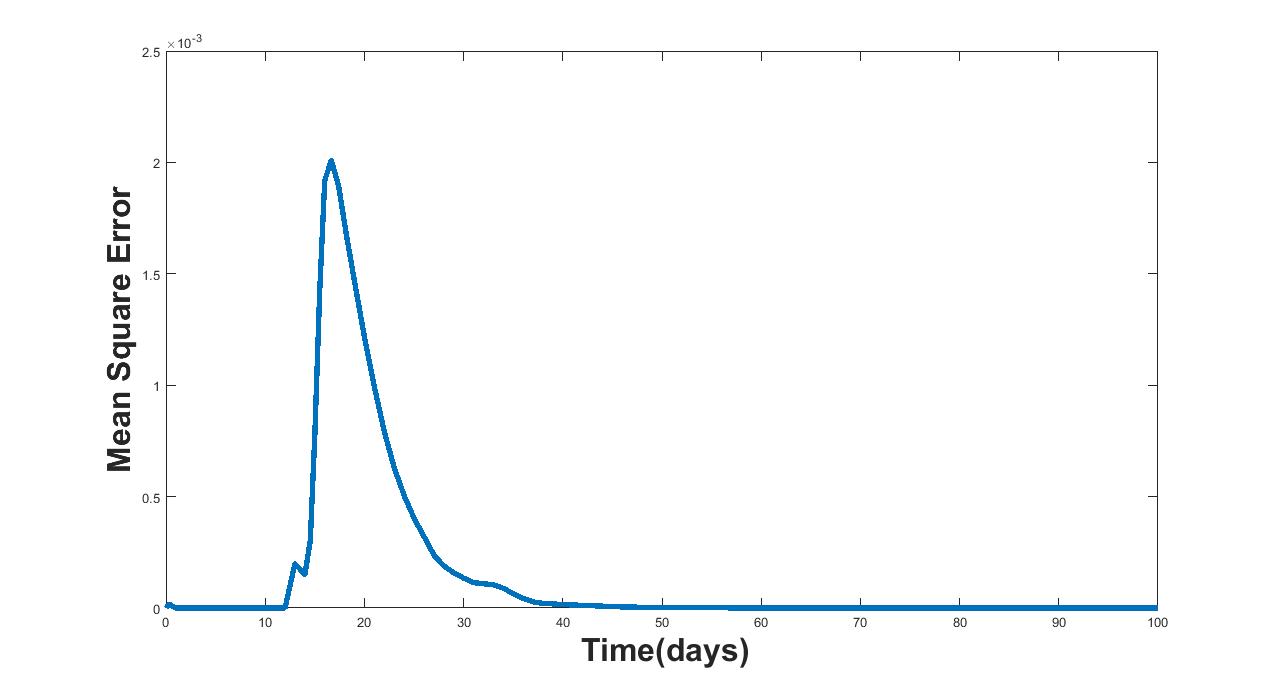}
\caption*{(a) Interval I : 0 to 1}
\end{center}
\end{figure}

\vspace{-3mm}
		
\begin{figure}[hbt!]
\begin{center}
\includegraphics[width=2in, height=1.8in, angle=0]{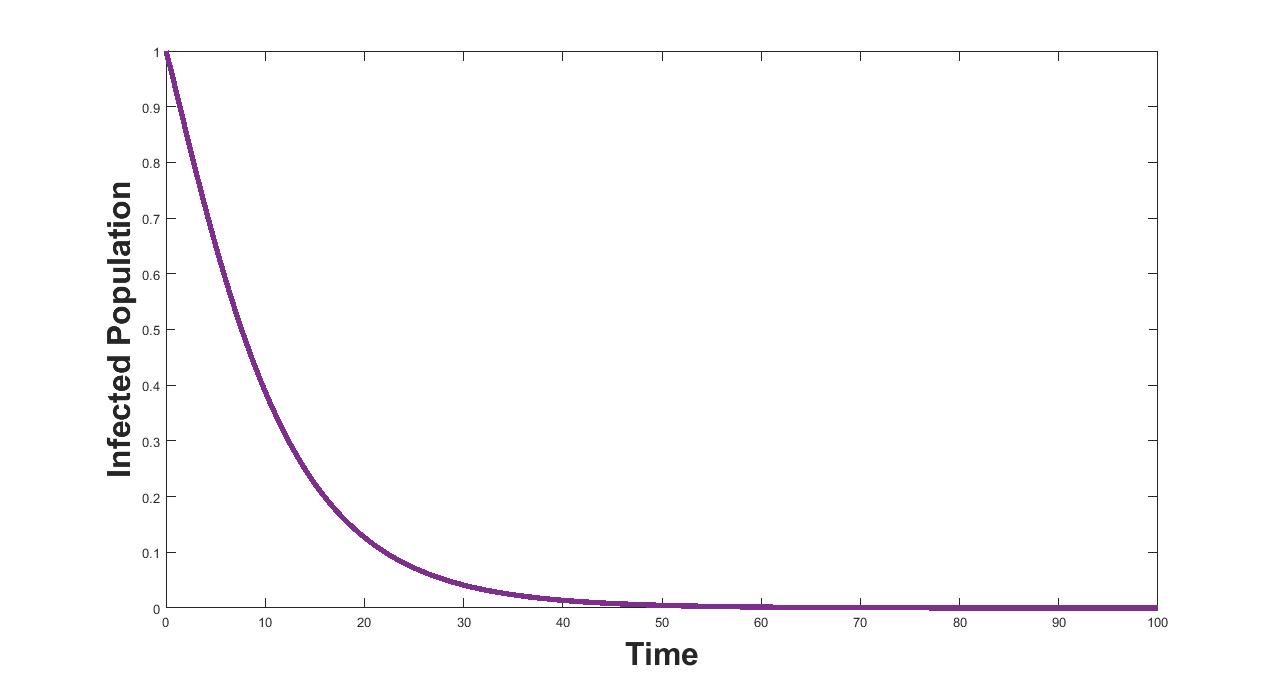}
\hspace{-.4cm}
\includegraphics[width=2in, height=1.8in, angle=0]{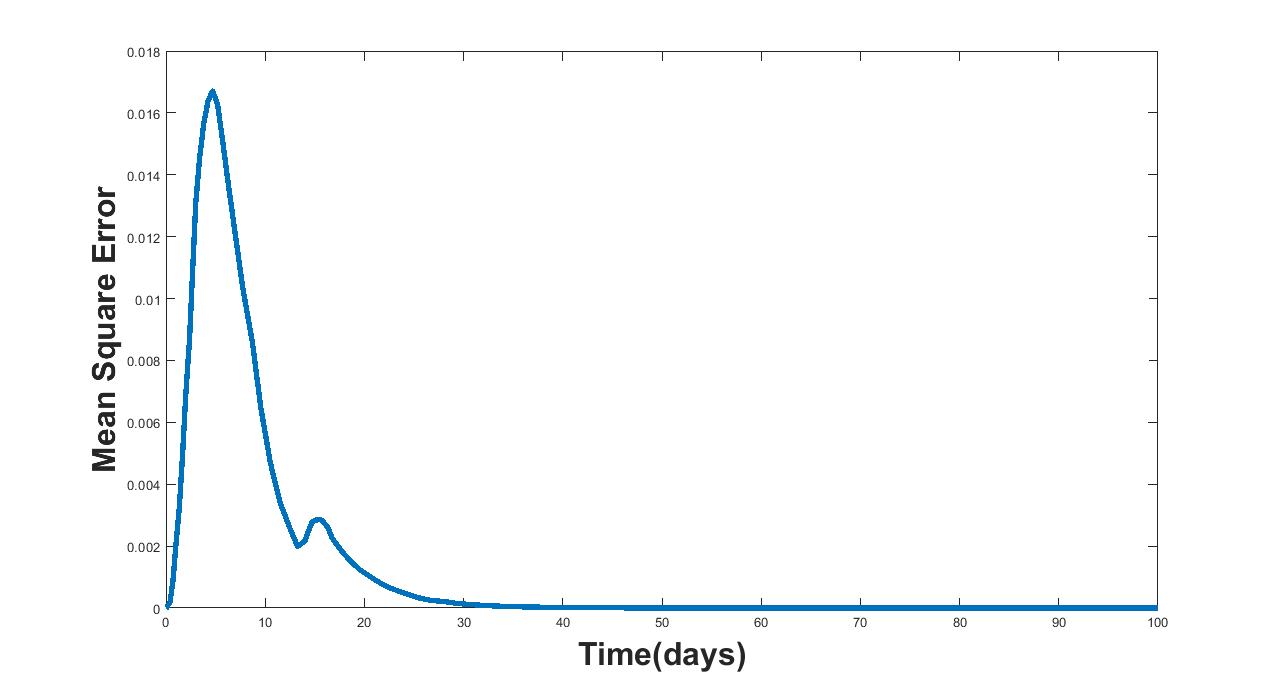}
\hspace{-.395cm}
\includegraphics[width=2in, height=1.8in, angle=0]{deltaerror2.jpg}
\caption*{(b) Interval I : 1 to 2.5}

\vspace{5mm}
\caption{Sensitivity Analysis of $\delta$}
\label{sen_delta}
\end{center}
\end{figure}

\newpage

\subsection{Parameter $\boldsymbol{\gamma}$}
	
\begin{figure}[hbt!]
\begin{center}
\includegraphics[width=2in, height=1.8in, angle=0]{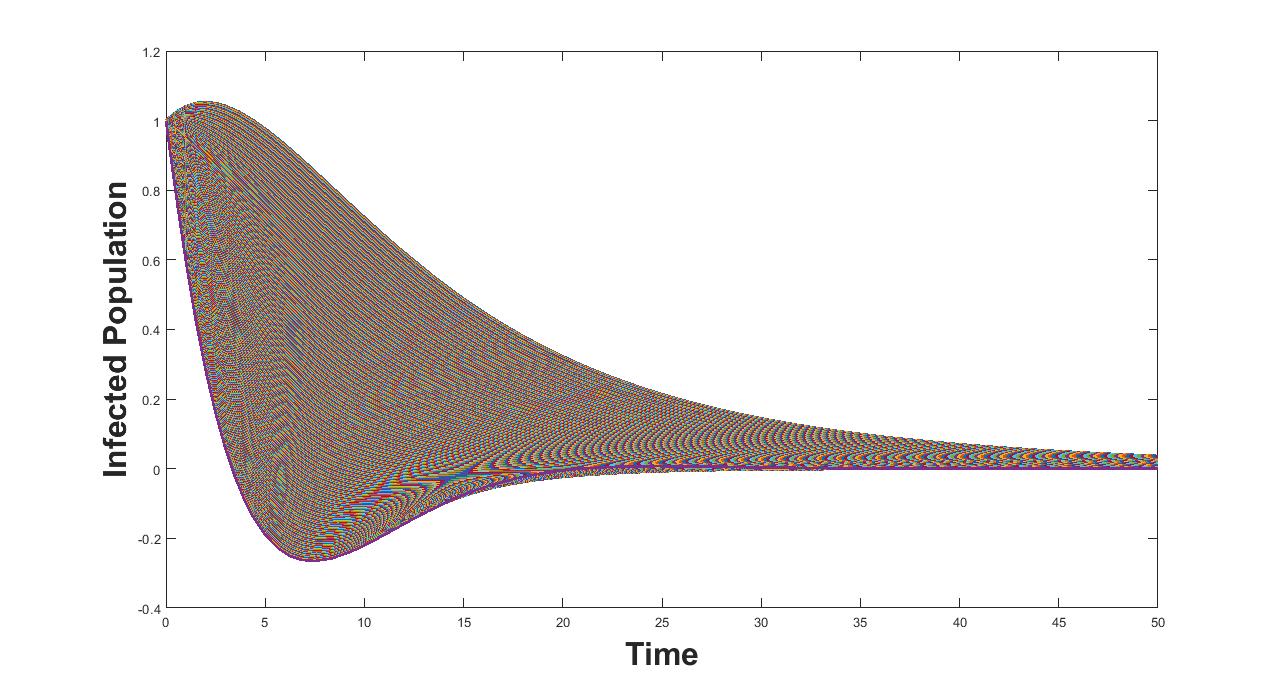}
\hspace{-.4cm}
\includegraphics[width=2in, height=1.8in, angle=0]{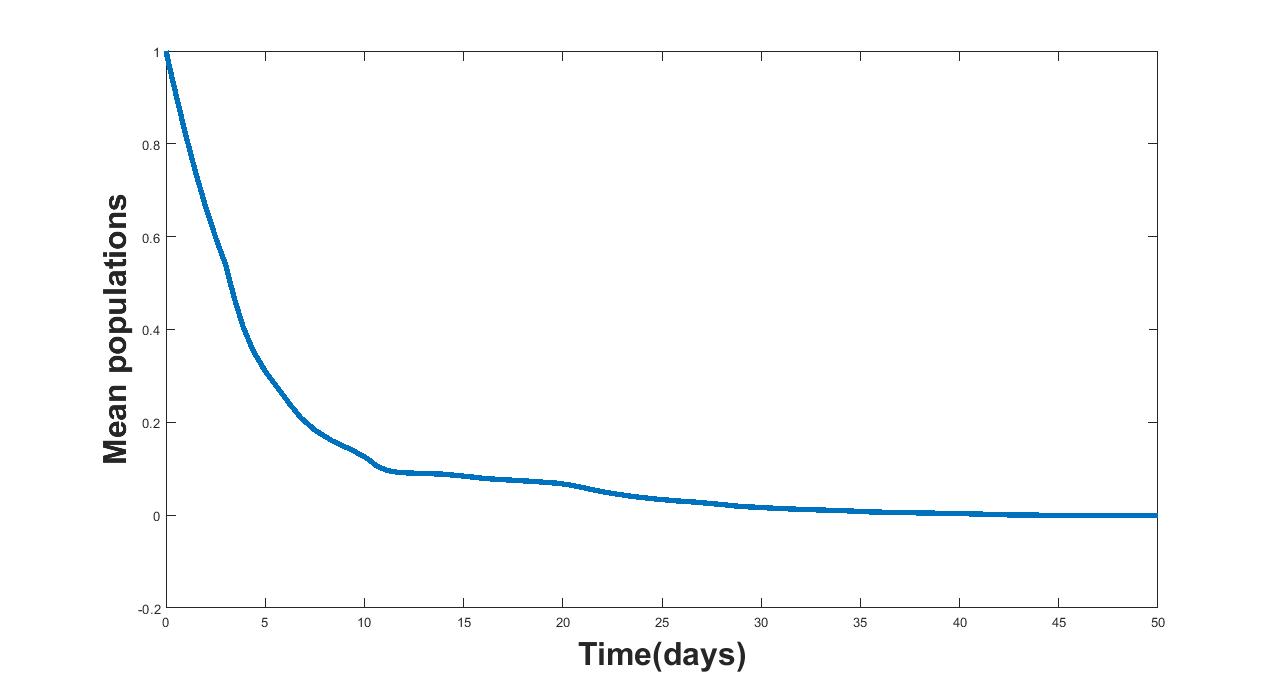}
\hspace{-.395cm}
\includegraphics[width=2in, height=1.8in, angle=0]{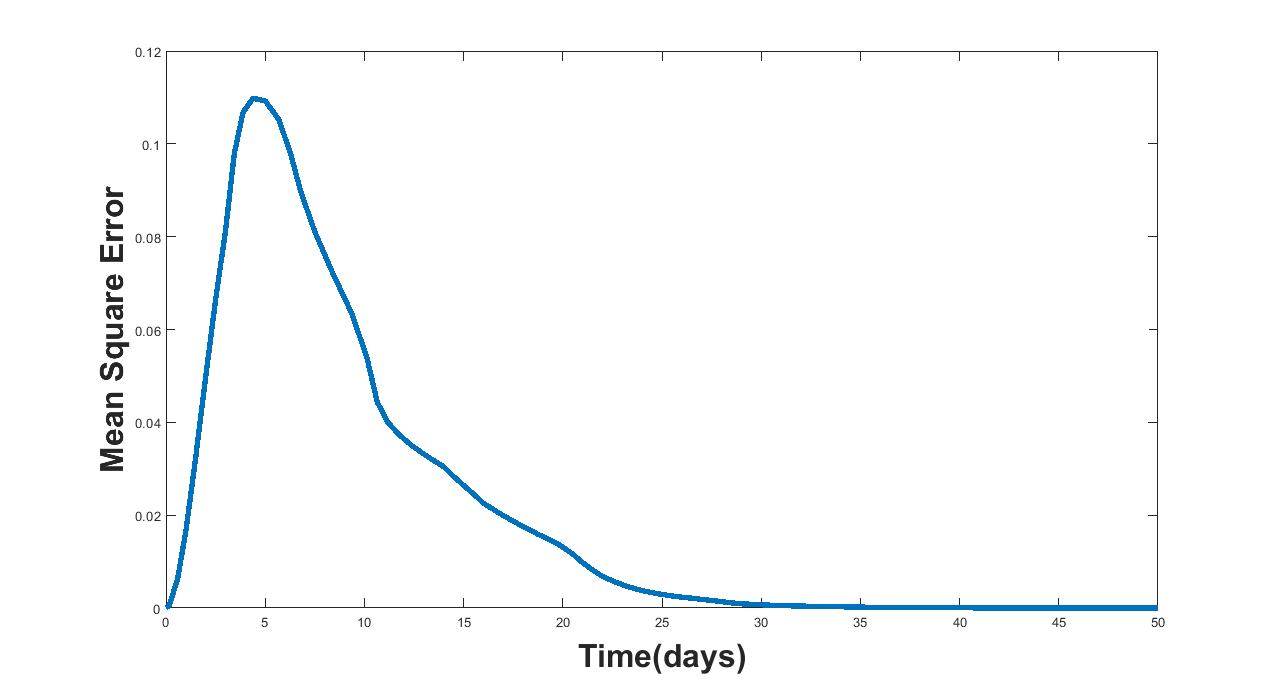}
\caption*{(a) Interval I : 0 to 1}
\end{center}
\end{figure}

\vspace{-3mm}
		
\begin{figure}[hbt!]
\begin{center}
\includegraphics[width=2in, height=1.8in, angle=0]{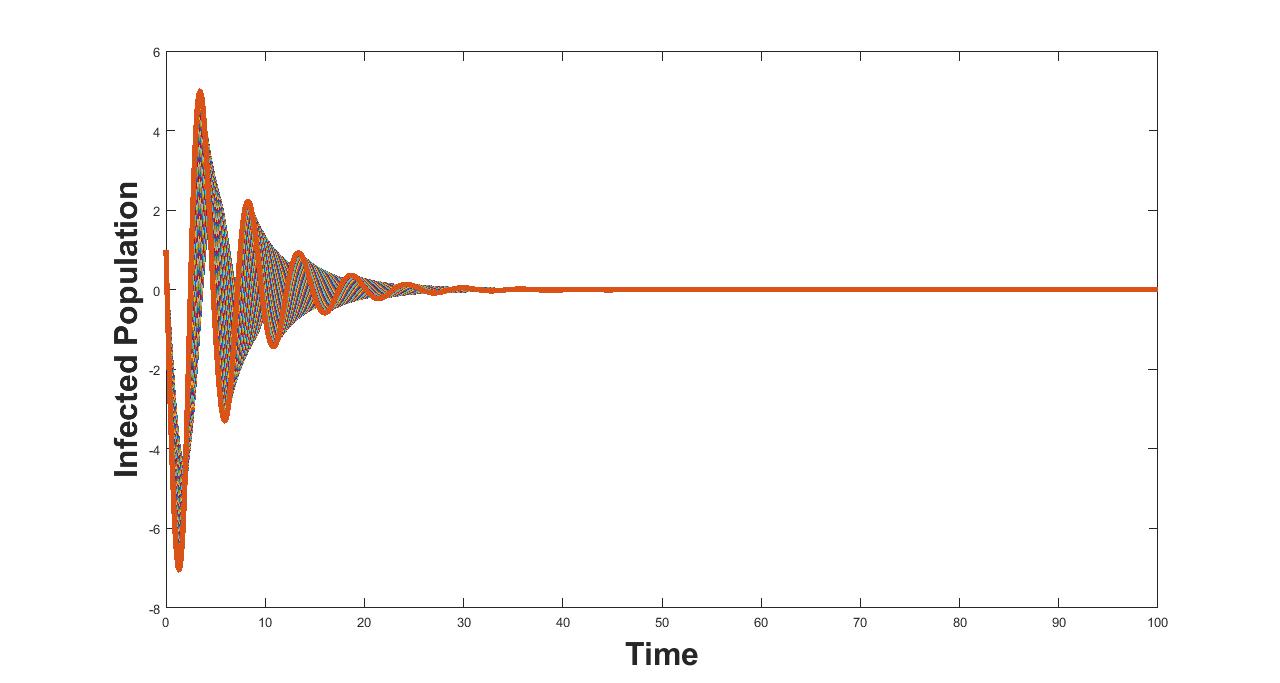}
\hspace{-.4cm}
\includegraphics[width=2in, height=1.8in, angle=0]{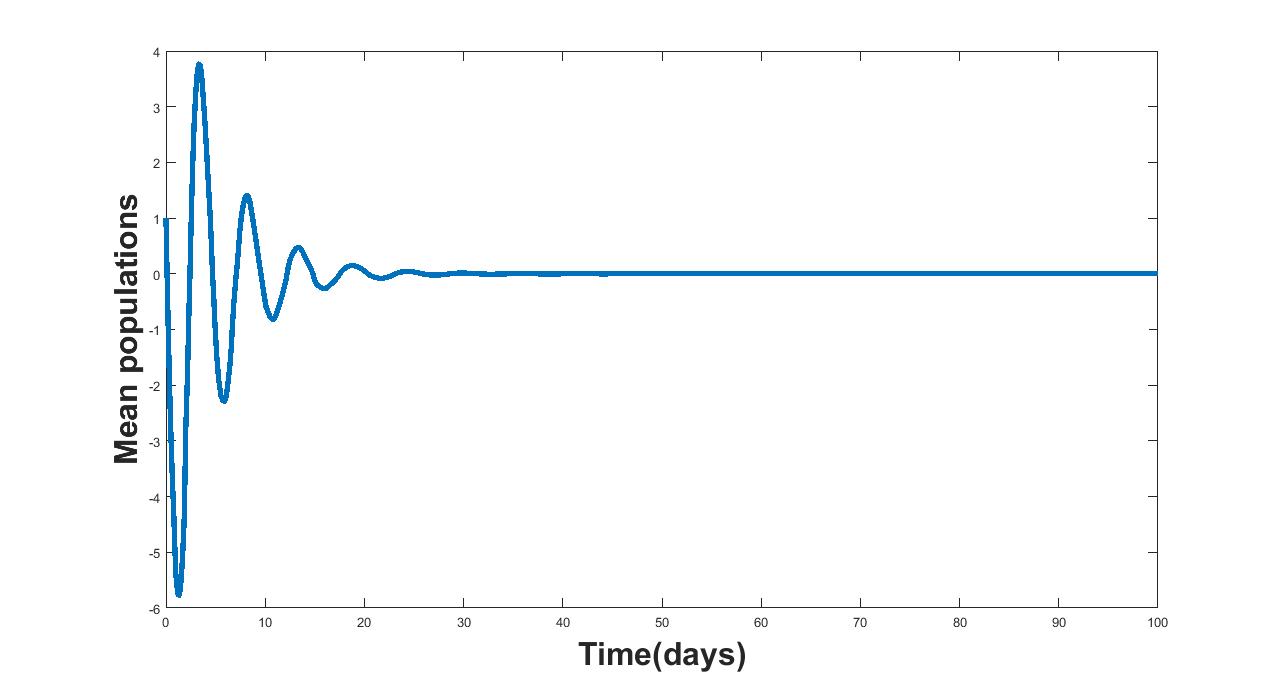}
\hspace{-.395cm}
\includegraphics[width=2in, height=1.8in, angle=0]{gammaerrorsen.jpg}
\caption*{(a) Interval I : 1 to 2.5}

\vspace{5mm}
\caption{Sensitivity Analysis of $\gamma$}
\label{sen_gamma}
\end{center}
\end{figure}

\end{document}